\documentclass[a4paper,12pt]{article}
\usepackage{amssymb}
\usepackage{amsbsy}
\usepackage{amsmath}
\usepackage{graphicx}
\newtheorem{theorem}{\bf Theorem}

\newtheorem{conjecture}[subsubsection]{\bf Conjecture}

\newtheorem{definition}[subsubsection]{\bf Definition}
\newtheorem{ex}[subsection]{\bf Example}
\newtheorem{proposition}[subsubsection]{\bf Proposition}
\newtheorem{lemma}[subsubsection]{\bf Lemma}
\newtheorem{corollary}[subsubsection]{\bf Corollary}
\newtheorem{rk}[subsubsection]{\bf Remark}
\newtheorem{rks}[subsubsection]{\bf Remarks}

\newenvironment{remark}{\begin{rk}\rm}{\end{rk}}

\newenvironment{pf}{{\em Proof. }}{\qed \newline}
\newenvironment{pfThm1i}{{\noindent \em Proof of Theorem
    \ref{new-mixing-theorem}(i).} }{\qed \newline}

\newenvironment{pf-asymp}{{\noindent \em Proof of Lemma
    \ref{asymp}.} }{\qed \newline}
\newenvironment{pfThm2-2}{{\noindent \em Proof of Theorem
    \ref{thm:perm-mix}(ii).} }{\qed \newline}

\newcommand{\BB}{\mathbb{B}}

\newcommand{\Z}{\mathbb{Z}}
\newcommand{\C}{\mathbb{C}}
\newcommand{\R}{\mathbb{R}}

\newcommand{\ee}{\mathbf{e}}
\newcommand{\vv}{\mathbf{v}}

\newcommand{\xx}{\mathbf{x}}
\newcommand{\yy}{\mathbf{y}}
\newcommand{\zz}{\mathbf{z}}

\newcommand{\qed}{\hspace*{\fill} \rule{2mm}{3mm}}

\newcommand{\lra}{\longrightarrow}
\newcommand{\LRA}{\Leftrightarrow}

\newcommand{\taumax}{\tau_\mathrm{max}}

\begin{document}

\title{On the mixing properties of piecewise expanding maps under composition with permutations}

\author{Nigel P.~Byott, Mark Holland, Yiwei Zhang, \\
College of Engineering, Mathematics and Physical Sciences, \\
University of Exeter, Exeter EX4 1QF, U.K.}

\date{\today}
\maketitle
\begin{abstract}

We consider the effect on the mixing properties of a piecewise
smooth interval map $f$ when its domain is divided into $N$ equal
subintervals and $f$ is composed with a permutation of these. The
case of the stretch-and-fold map $f(x)=mx \bmod 1$ for integers $m
\geq 2$ is examined in detail. We give a combinatorial description
of those permutations $\sigma$ for which $\sigma \circ f$ is still
(topologically) mixing, and show that the proportion of such
permutations tends to $1$ as $N \to \infty$. We then investigate the
mixing rate of $\sigma \circ f$ (as measured by the modulus of the
second largest eigenvalue of the transfer operator). In contrast to
the situation for continuous time diffusive systems, we show that
composition with a permutation cannot improve the mixing rate of
$f$, but typically makes it worse. Under some mild assumptions on
$m$ and $N$, we obtain a precise value for the worst mixing rate as
$\sigma$ ranges through all permutations; this can be made
arbitrarily close to $1$ as $N \to \infty$ (with $m$ fixed). We
illustrate the geometric distribution of the second largest
eigenvalues in the complex plane for small $m$ and $N$, and propose
a conjecture concerning their location in general. Finally, we give
examples of other interval maps $f$ for which composition with
permutations produces different behaviour than that obtained from
the stretch-and-fold map.
\end{abstract}

\section{Introduction}
Mixing processes of various kinds occur throughout nature and are
vital in many technological applications. It is therefore an
important and interesting problem to understand the properties of
these processes from a mathematical perspective.  In the context of
discrete time dynamical systems, transfer operator methods provide a
means of investigating such questions, and this approach has been
developed in a variety of settings \cite{Baladi94, Baladi00, BT05,
Boy97, GL06, HK89, Keller84, Viana97}. The transfer operator
$\mathcal{L}$ acts on a suitable Banach space of real-valued
functions (or distributions), and its spectrum provides a powerful
tool for analysing many mixing properties of the system, e.g.\
whether or not the system is indeed mixing \cite{Boy97}, the mixing
rate of the system \cite{Baladi00, Keller84, Viana97}, and the
existence of almost-invariant sets \cite{DJ99, FD}.

We restrict our attention to a piecewise smooth map $f$ on a compact
interval $I$.  We briefly recall
some facts about this situation; for more details see, for instance,
\cite{DFS00}. We consider the transfer operator
$\mathcal{L}_f$ of $f$ restricted to the Banach space $BV$ of
functions of bounded variation. The spectrum of
$\mathcal{L}_f|_{BV}$ is contained in the unit disk in the complex
plane. If $f$ is piecewise expanding, with the expansion factor
uniformly bounded away from 1, then the essential spectral radius
$r_{ess}$ can be interpreted as the slowest {\em local}
mixing rate of the system. The spectrum certainly contains the
eigenvalue $1$, corresponding to the equilibrium state of the system,
but may also contain further isolated points of modulus
greater than $r_{ess}$. These isolated eigenvalues come from
resonances in the system, and the corresponding eigenfunctions will converge to
equilibrium at a slower rate than would be predicted from $r_{ess}$.
The {\em global} mixing behaviour of the
system is therefore determined by the quantity
$\sup\left\{|\lambda|:\lambda\in\operatorname{Spec}(\mathcal{L}_{f}|_{\operatorname{BV}})\backslash\{1\}\right\}$. For
brevity, we will refer to this quantity as the {\em mixing rate} of the
system. Thus good mixing is achieved when the mixing rate is small,
while a mixing rate close to $1$ indicates that there are eigenfunctions
for which converge to equilibrium is very slow.

Even in the case of piecewise expanding interval maps, there seems to
be no general technique known for calculating the isolated eigenvalues
and hence finding the mixing rate.  A number of examples have
nevertheless been investigated in detail.  Baladi \cite{Baladi89}
constructed an expanding Markov map of constant slope for which the
transfer operator has a complex-conjugate pair of isolated
eigenvalues, and this was used by Collet and Eckman \cite{CE04} to
construct a two dimensional piecewise hyperbolic function which is the
(skew-)product of piecewise expanding interval map with similar
behaviour to Baladi's map. Dellnitz et al.\ \cite{DFS00} described a
parameterized family of expanding interval maps for which the location
of a non-trivial real, positive isolated eigenvalue may be controlled.

In this paper we study the effect on mixing of dividing $I$ into $N$
equal subintervals and composing $f$ with a permutation of these.
As far as we are aware, this is the first attempt to investigate the
effect of permutations on mixing for discrete dynamical systems. In
the continuous setting, Ashwin et al.\ \cite{Ashwin02} considered a
1-dimensional diffusion process, and showed that the mixing rate is
typically improved if subintervals of its domain are permuted at
regular time-steps. They considered permutations of various simple
kinds, and investigated numerically the effect of certain
permutations for small $N$. As well as treating the discrete
setting, the novelty of our approach is that we use combinatorial
and group-theoretic arguments to treat all permutations
systematically for arbitrarily large $N$.

We now describe more precisely the situation we investigate and the
main results we obtain. Let $f$ again be a piecewise smooth map on
the compact interval $I$, let $I$ be divided into $N$ subintervals
of equal length, and let $\sigma : I \to I$ be a piecewise smooth
map which simply permutes these intervals. (Thus we may identify
$\sigma$ with an element of the symmetric group $S_N$, which
consists of the $N!$ permutations of $N$ objects.) Then the
composite function $\sigma \circ f$ is again a piecewise smooth map
on $I$, and we wish to compare the global mixing behaviour of
$\sigma \circ f$ and $f$. The main focus of our study will be the
stretch-and-fold map $f(x)=mx \bmod 1$ on the interval $I=[0,1]$,
where $m \geq 2$ is an integer. This map is a standard (very simple)
example of a piecewise expanding interval map, and is itself often
taken as the canonical mixing protocol for polymers and pastes
\cite{Lasota94}. It can also be regarded as the prototype for the
much-studied family of maps $x \mapsto \beta x + \alpha \bmod 1$
(see for example \cite{Flatto97,Glendinning89, Hofbauer81}). Our
functions $\sigma \circ f$ provide a generalisation of the basic map
$f$ in a different direction.

For each choice of the two integer parameters $m$, $N$, we are
interested in the mixing behaviour of the collection of maps $\sigma
\circ f$ as $\sigma$ ranges through $S_N$. The mixing behaviour of $f$
itself is easy to describe. The essential spectral radius of
$\mathcal{L}_f$ is $1/m$, and there are no isolated eigenvalues
$\lambda$ with $|\lambda|> 1/m$ apart the simple eigenvalue 1.
Thus the mixing rate for $f$ is
$1/m$. To investigate the mixing behaviour of $\sigma \circ f$, we
must first address the issue of whether $\sigma \circ f$ is indeed
mixing at all. This is essentially a combinatorial question, depending
on $m$ and $N$ as well as on the particular permutation $\sigma$. We
will show in Theorem \ref{new-mixing-theorem} that, if $m$ is fixed,
then for many values of $N$, the function $\sigma \circ
f$ is mixing for {\em all} permutations $\sigma \in S_N$; for
the remaining values of $N$, the map $\sigma \circ f$ will fail to be
mixing for some permutations $\sigma$, but the proportion of such
permutations tends to $0$ as $N \to \infty$.

We will see that the essential spectral radius of
$\mathcal{L}_{\sigma \circ f}$ is again $1/m$, so that, when $\sigma
\circ f$ is mixing, its mixing rate can be no better than that of
$f$. For simplicity, we assume that $N>m$ and $\gcd(m,N)=1$. In particular,
this guarantees that $\sigma \circ f$ is mixing for all $\sigma \in
S_N$. Its mixing rate is
$$ \tau_\sigma :=
\sup\left\{|\lambda|:\lambda\in\operatorname{Spec}(\mathcal{L}_{\sigma\circ
f}|_{\operatorname{BV}})\backslash\{1\}\right\} \geq 1/m, $$ and
composition with $\sigma$ results in a {\em worse} mixing rate than
for $f$ alone unless we have equality. We determine in Theorem
\ref{thm:perm-mix} how bad the mixing rate can become: the maximal
value of $\tau_\sigma$ as $\sigma$ ranges through $S_N$ is
$\sin(m\pi/N) / m\sin(\pi/N)$, which can be made arbitrarily close
to $1$ by taking $N$ sufficiently large. The function $\sigma \circ
f$ is a Markov map, and an argument using Fredholm determinants
\cite{Mori89,Mori98} shows that $\tau_\sigma$ is the modulus of the
second largest eigenvalue of its probability transition matrix,
which is doubly stochastic. It is this which enables us to prove
Theorem \ref{thm:perm-mix}, the maximal value of $\tau_\sigma$ being
obtained when the matrix is conjugate to a circulant matrix.
Eigenvalues of various classes of stochastic matrices are
long-standing problems and have been discussed by many authors (see
for instance \cite{Berkolaiko01,
Dokovic,Ito,Karpelevich,Zyczkowski}). The proof of Theorem
\ref{thm:perm-mix} requires a result (Lemma \ref{circulant-worst})
on the effect of permuting the columns of a stochastic matrix; this
seems to be new and may be of independent interest. A natural
question is how, as $\sigma $ varies, the second largest isolated
eigenvalues of $\mathcal{L}_{\sigma\circ f}$ (and not just their
moduli) are distributed in the complex plane.  We propose a
conjecture on their distribution, on the basis some numerical
investigations.

The results just described relate to the particular maps $f(x)=mx
\bmod 1$, which are amenable to detailed combinatorial analysis. We
also briefly discuss two further cases, which exhibit different types
of behaviour. First, we exhibit a Markov map $f$ which is mixing, but
where the proportion of permutations $\sigma \in S_N$ with $\sigma
\circ f$ mixing does {\em not} tend to 1 as $N \to \infty$.  Secondly,
we give an example to show that, for a non-uniformly expanding map $f$
with intermittent behaviour, composition with permutations may speed
up the mixing rate.  These two examples indicate that one cannot
expect results along the lines of our Theorems
\ref{new-mixing-theorem} and \ref{thm:perm-mix} to hold for arbitrary
piecewise expanding interval maps $f$. Nevertheless, the results we
have obtained suggest that the effect of composition with permutations
is fundamentally different for discrete and continuous dynamical
systems: it typically results in improved mixing in the continuous
case, but worse mixing in the discrete case.

The organisation of this paper is as follows. In \S\ref{background},
we give the necessary background and then state our main results. We
also briefly discuss the location of the isolated eigenvalues in the
complex plane. The proof of Theorem \ref{new-mixing-theorem} is given
in \S\ref{sec-t1}, along with some explicit formulae for the
proportion of non-mixing permutations in special cases. This section
is essentially combinatorial in character. Theorem \ref{thm:perm-mix}
is proved in \S\ref{sec_proof}. Finally, the two additional examples
mentioned above are presented in \S\ref{examples}.

\section{Background and statement of results}  \label{background}
\subsection{Mixing versus non-mixing}\label{sec_mixvsnonmix}
In this section we state our main result in relation to the question of mixing versus non-mixing of
$\sigma\circ f$. Given a measure preserving system $(f,M,\mu)$, we say that the system is (strongly) mixing if
\begin{equation}\label{mixing_def}
|\mu(f^{-n}A\cap B)-\mu(A)\mu(B)|\to 0 \mbox{ as } n\to\infty,
\end{equation}
where $A,B$ are $\mu$-measurable sets. Another version of mixing is
that of \emph{topological mixing}, namely we say $(f,M)$ is
topologically mixing if for all open $U,V\subset M$, there exists a
constant $n_0=n_0(U,V)$ such that $\forall~n\geq n_0$, $f^n(U)\cap
V\neq\emptyset$. To show that $f$ is \emph{not} mixing, it is
usually easier to show that $f$ is not
topologically mixing.

For the examples that we consider, it will be also true that
topological mixing implies strong mixing, see \cite{Viana97}.

We will consider maps on the unit interval $I$, dividing $I$ into $N$
equal subintervals. To avoid the problem of functions
being undefined, or multiply defined, at endpoints of these
subintervals, we work with (non-compact) intervals which are closed on
the left and open on the right. Thus we consider piecewise continuous
maps $f \colon [0,1) \lra [0,1)$. We can of course regard $f$ as a map
    on the compact interval $[0,1]$ (by stipulating $f(1)=f(0)$) or on
    the circle $\mathbb{S}^1=\R/\Z$.

We divide the unit interval as follows. Fix $N \geq 2$, and let $I_j=[j/N, (j+1)/N )$, $0 \leq j
    <N$. For any permutation $\sigma$ of $\{0, 1, \ldots, N-1\}$ we
    write $\sigma$ also for the corresponding interval exchange map:
$$ \sigma(x) = x + (\sigma(j)-j)/N \bmod 1 \mbox{ for } x \in I_j. $$
We write $S_N$ for the group of all permutations of $\Z/N\Z$.

The specific map $f \colon [0,1) \lra [0,1)$ we consider is $f(x)=mx
\bmod 1$ for a fixed integer $m \geq 2$.  Our first result shows
that the composite $\sigma \circ f$ is mixing for {\em almost all}
permutations $\sigma$ when $N$ is large enough.

\begin{theorem} \label{new-mixing-theorem}
Let $f$ be as above. Then
\begin{enumerate}
\item[(i)] if $N$ is not a multiple of $m$ then
$\sigma \circ f$ is mixing for all $\sigma \in S_N$;
\item[(ii)] if $N>m$ and $N$ is a multiple of $m$, say $N=m \ell$, then there will be
  some $\sigma \in S_N$ for which $\sigma \circ f$ is not mixing. As
  $\ell \to \infty$ (with $m$ fixed), however, the
  proportion of permutations $\sigma$ with $\sigma \circ f$ mixing
  tends to $1$.
\end{enumerate}
\end{theorem}

The proof of Theorem \ref{new-mixing-theorem} is given in
\S\ref{sec-t1}.



\subsection{Background on transfer operators}\label{backgroud}
Our methods for studying mixing rates will utilise the theory of
transfer operators and Fredholm matrices, see \cite{Baladi00,
Keller84, Mori89, Mori06, Viana97}. We now give an overview of the
relevant theory. For a measure preserving system $(f,M,\mu)$, the
rate of mixing can be quantified in various ways. However, we will
primarily focus on the \emph{speed of convergence to equilibrium}.
More precisely, if $f:[0,1]\to[0,1]$ is a piecewise expanding map,
we define the transfer operator $\mathcal{L}_{f}:L^1\to L^1$ by:
\begin{equation}
\{\mathcal{L}_{f}\phi\}(x)=\sum_{f(y)=x}\frac{\phi(y)}{|f'(y)|},\quad\forall\phi\in
L^1.
\end{equation}
The operator $\mathcal{L}_f$ satisfies the following identity, for
$\phi\in L^p,\psi\in L^q$ (with $p^{-1}+q^{-1}=1$):
\begin{equation}
\int(\mathcal{L}_{f}\phi)\psi\,dx=\int\phi(\psi\circ f)\,dx,
\end{equation}
where $dx$ denotes integration with respect to the reference
(Lebesgue) measure. If $f$ preserves an ergodic measure $\mu$ with
density $\rho(x)\in L^1$, then $(\mathcal{L}_f\rho)(x)=\rho(x)$.
Suppose now that we have a Banach space $\mathcal{B}\subset L^1$,
with $\rho\in\mathcal{B}$, and with norm $\|\cdot\|_{\mathcal{B}}$.
We define the speed of convergence to equilibrium in $\mathcal{B}$
as the rate $r(n)$ such that $\exists~C_{\mathcal{B}}<\infty$,
\begin{equation}\label{eq:rate-rn1}
\|\mathcal{L}^{n}_{f}(\phi)-\rho\|_{\mathcal{B}}\leq
C_{\mathcal{B}}r(n),\quad\forall\phi\in\mathcal{B},\,\|\phi\|_1=1,\,\forall
n\geq 1,
\end{equation}
and there exists $\phi\in\mathcal{B}$ with $\|\phi\|_1=1$, such that,
for sufficiently large $n$ and for some $C_\phi>0$, we have
\begin{equation}\label{eq:rate-rn2}
\|\mathcal{L}^{n}_{f}(\phi)-\rho\|_{\mathcal{B}}\geq
C_{\phi}r(n),\quad\forall\phi\in\mathcal{B}.
\end{equation}
For the whole space $L^1$, the rate function $r(n)$ cannot be
specified, i.e. there exist $\phi\in L^1$ for which
$\|\mathcal{L}^{n}_{f}(\phi)-\rho\|_1$ decays arbitrarily slowly.
When $f$ is a piecewise expanding map, the natural space to consider
is $\mathcal{B}=\operatorname{BV},$ the class of functions with
\emph{bounded variation}. We recall this definition as follows, see
\cite{Keller84}. Given a function $\phi:[0,1]\to\mathbb{R}$, we
define the total variation of $\phi$ as
\begin{equation}
\textrm{var}(\phi)=\sup\{\sum_{k=1}^{n}|\phi(x_k)-\phi(x_{k-1})|:0\leq x_0\leq\ldots\leq x_n=1\},
\end{equation}
where the $\sup$ is taken over all partitions of $[0,1]$. We say
that $\phi$ has bounded variation (i.e. $\phi\in \operatorname{BV}$)
if $\textrm{var}(\phi)<\infty$. To make $\operatorname{BV}$ into a
Banach space, we define the norm $\|\cdot\|_{\operatorname{BV}}$ by
\[
||\phi||_{\operatorname{BV}}:=||\phi||_{1}+\operatorname{var}(\phi),
\]
and hence consider functions $\phi\in \operatorname{BV}$ with
$\|\phi\|_{\operatorname{BV}}<\infty$. Bounds on the rate of mixing
$r(n)$ can then be obtained by analysing the spectral properties of
the restriction $\mathcal{L}_{f}|_{\operatorname{BV}}:
\operatorname{BV}\to \operatorname{BV}$. In particular we say that
$\mathcal{L}_{f}|_{\operatorname{BV}}$ has a {\em spectral gap} if
\begin{equation}\label{equ_mixingrate}
\tau:=\sup\left\{|\lambda|:\lambda\in\operatorname{Spec}(\mathcal{L}_{f}|_{\operatorname{BV}})\backslash\{1\}\right\}<1,
\end{equation}
where $\operatorname{Spec}(\mathcal{L}_{f})$ is the spectrum of $\mathcal{L}_f$.
Hence for any $\phi\in\operatorname{BV}$ and any $\epsilon>0$, the
spectral decomposition of $\mathcal{L}_{f}|_{\operatorname{BV}}$
implies that
there exists $C>0$ such that for all $n$,
\begin{equation}\label{equ_spectraldeco}
||\mathcal{L}^{n}_{f}(\phi)-\rho||_{\operatorname{BV}}\leq
C\cdot(\tau+\epsilon)^{n}||\phi||_{\operatorname{BV}}\;\textrm{with}\;\|\phi\|_1=1.
\end{equation}
Thus $\tau$ determines the rate of convergence to
equilibrium. As noted in the Introduction, we will refer to $\tau$ as
the  \emph{mixing rate} of $f$.
Since there is no general method for finding the exact value of
$\tau$, we consider also the \emph{essential spectral
radius}, $r_{ess}=r_{ess}(\mathcal{L}_{f}|_{\operatorname{BV}})$, defined by
\[
r_{ess}
:=\inf\{r\geq 0:
\lambda\in\operatorname{Spec}(\mathcal{L}_{f}|_{\operatorname{BV}}),\,|\lambda|>r\implies\lambda\,\textrm{isolated}\}.
\]
The isolated eigenvalues $\lambda$ with $|\lambda|>r_{ess}$
are of finite multiplicity. For piecewise expanding maps, see
\cite[Theorem 1]{Keller84},
we have the lower bound on $\tau$ via
\begin{equation}\label{equ_inequility}
\tau\geq
r_{ess}=\exp\left\{-\liminf_{k\to\infty}\operatorname{essinf}\limits_{x\in[0,1]}\frac{1}{k}\log|(f')^{k}(x)|\right\}.
\end{equation}

\subsection{The main result on mixing rates} \label{sec3:statement}
In our setting we consider specifically the map $f(x)=mx\mod 1$.
When $f$ is composed with a permutation $\sigma\in S_N$, we have seen
in Theorem \ref{new-mixing-theorem} that the
resulting piecewise linear transformation $\sigma \circ f$ is usually
(but not always) mixing. When $\sigma \circ f$ is mixing, we consider
its mixing rate
\begin{equation} \label{def-tau-sig}
\tau_\sigma :=
\sup\left\{|\lambda|:\lambda\in\operatorname{Spec}(\mathcal{L}_{\sigma
  \circ f}|_{\operatorname{BV}})\backslash\{1\}\right\}.
\end{equation}
We state the following result.

\begin{theorem}\label{thm:perm-mix}
Fix $m$, $N \geq 2$ and consider the transformations $\sigma \circ f$
where $f(x)=mx\mod 1$ and $\sigma \in S_N$. Then the following hold.
\begin{enumerate}
\item[(i)] For all $\sigma\in S_N$, the essential spectral radius is given by $r_{ess}(\mathcal{L}_{\sigma\circ f}|_{\operatorname{BV}})=1/m$.
\item[(ii)] If $N>m$ and $\gcd(m,N)=1$, then, for each $\sigma \in
  S_N$, we have
$$  \tau_\sigma \leq \taumax :=   \frac{\sin(m\pi/N)}{m\sin(\pi/N)} .  $$
Moreover, each of the values
$(-1)^{m-1} e^{2 \pi i j/N} \taumax$ for $0 \leq j <N$ and $(-1)^m
\taumax$ occurs as an isolated
  eigenvalue of $\mathcal{L}_{\sigma \circ f}$ for an appropriate
    choice of $\sigma$. Thus $\tau_\sigma=\taumax$ for these $\sigma$.
\end{enumerate}
\end{theorem}
The proof of Theorem \ref{thm:perm-mix}
is given in \S\ref{sec_proof}.

\subsection{Geometric location of the isolated eigenvalues}

For given $m$ and $N$, and for $\sigma \in S_N$, let
$$ \Lambda_{\sigma} :=
\{\lambda\in\operatorname{Spec}(\mathcal{L}_{\sigma\circ
  f}|_{\operatorname{BV}})~\mbox{such that}~|\lambda|=\tau_{\sigma}\}, $$
where $\tau_{\sigma}$ is defined in (\ref{def-tau-sig}). We
will see below that $\Lambda_\sigma$ is the set of second largest
eigenvalues of a certain doubly stochastic matrix, namely the
probability transition matrix for the Markov map $\sigma \circ
f$. (Here, ``second largest'' is in terms of the modulus, the largest
eigenvalue always being $1$, and a given matrix may have more than one
second largest eigenvalue since there may be distinct eigenvalues with
the same modulus.)

\begin{figure}[h]
  \centering
  \includegraphics[width=16cm]{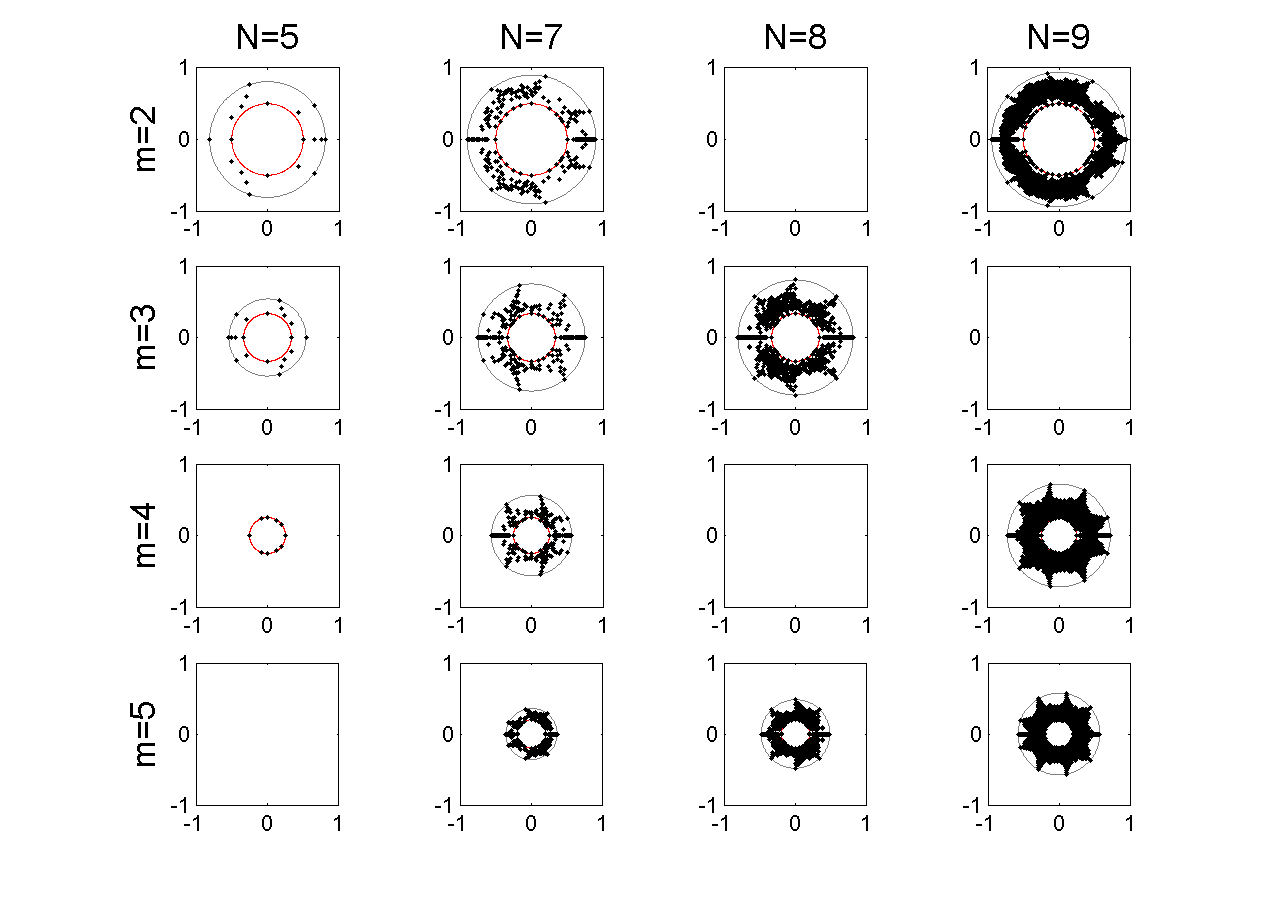}\\
  \caption{Geometric location of second largest isolated eigenvalues for the composition $\sigma\circ f$ where $f(x):=mx\mod 1$ and $\sigma\in S_{N}$ with $\operatorname{gcd}(N,m)=1$.}
  \label{fig_eigenvaluepuzzle}
\end{figure}

We would like to understand the geometric properties of the (finite)
set $\bigcup_{\sigma \in S_N} \Lambda_\sigma$ in the complex plane.
The elements of this set are the isolated eigenvalues of
$\mathcal{L}_{\sigma\circ f}$ which determine the mixing rates of
the maps $\sigma \circ f$ for all permutations $\sigma \in S_N$. For
small values of $m$ and $N$ with $\gcd(m,N)=1$, these sets are shown
in Figure \ref{fig_eigenvaluepuzzle}. In particular, these
eigenvalues are located between the inner circle with radius $1/m$
and the outer circle with radius $\taumax$.  This is in agreement
with Theorem \ref{thm:perm-mix}. As a special case, when $m=N-1=4,$
so $\sin(m\pi/N)=\sin(\pi/N)$, the two circles coincide.

The set $\mathcal{D}_N$ of doubly stochastic matrices of order $N$ has
good convexity properties \cite[Chapter I, \S5]{Schaefer74}. By a well-known
result of Birkhoff, $\mathcal{D}_N$ is precisely the convex hull of
the permutation matrices of order $N$. Moreover, the eigenvalues of
matrices in $\mathcal{D}_N$ lie in the convex hull of all roots of
unity of order at most $N$.


Together with Figure \ref{fig_eigenvaluepuzzle}, this suggests the
following conjecture:

\begin{conjecture}\label{conjecture_convex}
Suppose that $\gcd(m,N)=1$. Then $\bigcup_{\sigma\in
S_{N}}\Lambda_{\sigma}$ is contained in the convex hull of the
points $(-1)^{m-1} e^{2 \pi i j/N} \taumax$ for $0 \leq j <N$ and
$(-1)^m\taumax$. In particular, these are the only points $\lambda
\in \bigcup_{\sigma\in   S_{N}}\Lambda_{\sigma}$ with
$|\lambda|=\taumax$.
\end{conjecture}

We note that the convex hull in Conjecture \ref{conjecture_convex}
is a regular $N$-gon if $N$ is even, and an irregular $(N+1)$-gon
(obtained by adding one extra vertex to a regular $N$-gon) if $N$ is
odd; c.f.\ Figure \ref{fig_eigenvaluepuzzle}.

\section{Permutations preserving mixing for $mx \bmod 1$} \label{sec-t1}

Our main goal in this section is the proof of Theorem
\ref{new-mixing-theorem}. In \S\ref{sect-non-mixing}, we prove
statement (i), and, in the setting of statement (ii), give
a group theoretic interpretation of those $\sigma$ for which $\sigma
\circ f$ is not mixing. The asymptotic analysis needed to complete the
proof of Theorem \ref{new-mixing-theorem} is given in
\S\ref{l-asymp}.  In \S\ref{l-small} we give explicit formulae for the
proportion of non-mixing permutations for small values of $\ell$.

\subsection{When is $\sigma \circ f$ non-mixing?}
\label{sect-non-mixing}

Recall that
 $f(x) = mx \bmod{1}$ and $\sigma \in S_N$, and that we partition the
unit interval into subintervals $I_{a}=[a/N,(a+1)/N)$ for $a \in \{0,
  1, \ldots, N-1\}$.
We identify the indexing set $\{0, \ldots, N-1\}$ with the
the ring $\Z/N\Z$ of integers modulo $N$, so that arithmetic in this
indexing set is to be interpreted as arithmetic modulo $N$.

To begin with, we allow arbitrary $m$, $N \geq 2$. We set $g=\sigma
\circ f$.

\begin{definition} \label{sig-tilde}
For any subset $A \subseteq \Z/N\Z$, we define
$$ \tilde{f}(A) = \bigcup_{d=0}^{m-1} (mA+d) \subseteq  \Z / N\Z, \qquad
   \tilde{g}(A) = \sigma(\tilde{f}(A)). $$
\end{definition}

\begin{proposition} \label{g-union}
For each $A \subseteq \Z/N\Z$, we have
$$ f\left( \bigcup_{a \in A} I_a \right) =
      \bigcup_{b \in \tilde{f}(A)} I_b, \qquad
g\left( \bigcup_{a \in A} I_a \right) =
      \bigcup_{b \in \tilde{g}(A)} I_b. $$
\end{proposition}

\begin{pf}
This is immediate since for
each $j \in \Z/N\Z$ we have
$$ f(I_j) = \bigcup_{d=0}^{m-1} I_{mj+d}, \qquad \sigma(I_j) =
I_{\sigma(j)}. $$
\end{pf}

\begin{proposition} \label{size-At}
For each $A \subseteq \Z/N\Z$, we have
$$ \sharp A \leq \sharp \tilde{f}(A)  \leq m \sharp A. $$
Moreover, suppose that $0 < \sharp A  < N$. Then we have $\sharp
\tilde{f}(A)  =\sharp A$
if and only if the following two conditions hold:
\begin{itemize}
\item[(i)] $N=m \ell$ for some integer $\ell$;
\item[(ii)] $A$ is a union of cosets of $\ell\Z / N \Z$ (that is, $j
  \in A \Rightarrow j+\ell \in A$ for all $j \in \Z/N\Z$).
\end{itemize}
\end{proposition}
\begin{pf}
If $b \in \tilde{f}(A)$ then
$ b \equiv ma+d \pmod{N}$ for at least one of the $m \sharp A$ pairs $(a,d)$ with $a \in A$
and $0 \leq d <m$. Hence
$\sharp \tilde{f}(A)  \leq m \sharp A$. Now fix one pair
$(a_0,d_0)$. If another pair $(a,d)$ gives the same element $b$ then
\begin{equation} \label{congr}
           ma + d \equiv m a_0 + d_0 \pmod{N}.
\end{equation}
Thus $d \equiv d_0 \pmod{s}$, where $s=\gcd(m,N)$. This gives $m/s$
possibilities for $d$. For each of these, \eqref{congr} has $s$
solutions $a$ in $\Z/N\Z$, all
congruent mod $N/s$ (but in general not all in $A$).  So each $b$
arises from at most $m$ of the pairs
$(a,d)$, giving $\sharp \tilde{f}(A) \geq \sharp A$. This proves the first assertion.

If $\sharp \tilde{f}(A)= \sharp A$, then each $b$ must
arise from exactly $m$ pairs $(a,d)$. Thus given $a_0 \in
A$, we may take $d=d_0=0$, and the $s$ solutions $a$ to (\ref{congr})
must all lie in $A$. This shows that $a_0 + N/s \in A$, so that $A$ is
stable under addition of $N/s$.

First suppose (i) holds. Then $N/s=\ell$, so that if $\sharp \tilde{f}(A)= \sharp A$ then
(ii) holds. Conversely, if (i) and (ii) hold, then each $b \in \tilde{f}(A)$
arises from $m$ pairs $(a+j\ell,d)$ with $0 \leq j <m$, so that
$\sharp \tilde{f}(A)= \sharp A$.

It remains to show that if (i) does not hold and $\sharp
\tilde{f}(A)=\sharp A>0$
then $A=\Z / N\Z$. So let $m=es$ with $e>1$, and let
$a_0 \in A$. Since $s<m$, we may take $d_0=s$ in (\ref{congr}). But
(\ref{congr}) must have $s$ solutions for each of the possible values
$d \equiv d_0 \pmod{s}$ with $0 \leq d <m$, so we can find
$a_1 \in A$ with $ma_1 \equiv ma_0+s \pmod{N}$. Then $e a_1 \equiv e
a_0 +1 \pmod{N/s}$. Iterating, we can find $a_j \in A$ with $e a_j
\equiv ea_{j-1}+1 \equiv ea_0+j \pmod{N/s}$ for $j \geq 1$. As
$\gcd(e,N/s)=1$, we have $a_e \equiv a_0 + 1 \pmod{N/s}$. Since we
already know that $A$ is stable under addition of $N/s$, it follows
that $A$ is stable under addition of 1, so that $A=\Z/N\Z$.
\end{pf}

\begin{corollary} \label{g-til}
For any $A \subseteq \Z / N\Z$, we have
$$ \sharp \tilde{g}(A) \geq \sharp A.  $$
Moreover, if $A$ is a proper subset of $\Z/N \Z$ then equality can
only occur if $N=\ell m$ for some integer $\ell$.
\end{corollary}
\begin{pf}
This is clear since $\sharp \tilde{g}(A)=\sharp \tilde{f}(A)$.
\end{pf}

\begin{lemma}  \label{tilde-crit}
$g$ fails to be (topologically) mixing if and only if there is some proper subset $A$
  of $\Z/N\Z$ such that $\sharp \tilde{g}^r(A) = \sharp A$ for all $r
  \geq 0$.
\end{lemma}
\begin{pf}
Let $A$ be a subset with $0<\sharp A<N$ and $\sharp \tilde{g}^r(A) =
\sharp A$ for all $r$. As there are only finitely many subsets of $\Z/N\Z$,
we may choose $s \geq 0$ and $t \geq 1$ with
$\tilde{g}^{s+t}(A) =\tilde{g}^s(A)$.  Set
$B=\tilde{g}^s(A)$ and take non-empty open sets $U \subset
I_j$ and $V \subset I_k$ where $j \in B$ and $k \not \in B$. Then for
all $n \geq 0$ we have $g^{nt}(U) \subseteq \bigcup_{b \in B} I_b$
so that $g^{nt}(U) \cap V = \emptyset$. Hence $g$ is not mixing.

Conversely, suppose there is no proper subset $A$ with
$\sharp \tilde{g}^r(A) = \sharp A$ for all $r$. To see that $g$ is
mixing, we show that for any non-empty open subset $U$ of $[0,1)$ we have
 $g^n(U)=[0,1)$ for large enough $n$. Without loss of generality, $U$
is an interval of length $\delta>0$. Since $m>1$, we can choose $h$
large enough that $g^h(U)$ contains the initial point $j/N$ of some interval
$I_j$. Then for some $\epsilon>0$, we have $[j/N,j/N+\epsilon)
\subseteq g^h(U) \cap I_j$. Choose $k$ so that $m^k\epsilon >
1/N$ and let $g^k(j/N)=j'/N$. Then $I_{j'} \subset g^k(I_j) \subset
g^{h+k}(U)$. Now let $B=\{j'\}$ and take $s\geq 0$, $t \geq 1$ with
$\tilde{g}^{s+t}(B) =\tilde{g}^s(B)$.  The
non-empty set $A=\tilde{g}^s(B)$ then satisfies the
condition $\tilde{g}^{nt}(A) = A$ for all $t \geq 0$. Hence, by
Corollary \ref{g-til}, we have $\sharp \tilde{g}^r(A) =\sharp A$
for all $r \geq 0$. Thus our
hypothesis forces $A=\Z/N\Z$, so that $g^q(I_{j'}) = [0,1)$ for
all $q \geq s$. It follows that $g^n(U) = [0,1)$ for all $n \geq
h+k+s$, as required.
\end{pf}

\begin{pfThm1i}
Suppose that $N$ is not a multiple of $m$, and let $g=\sigma \circ f$
with $\sigma \in S_N$. By Corollary \ref{g-til} there is no proper subset $A$
with $\sharp \tilde{g}(A)  =\sharp A$. Hence by Lemma
\ref{tilde-crit}, $g$ is mixing.
\end{pfThm1i}

We now suppose that the $N=m\ell$ for some integer $\ell \geq 1$.

\begin{proposition} \label{delta-def}
There exists a permutation $\delta \in S_N$ such that
\begin{equation} \label{delta-prop}
 f(I_j) \supseteq I_{\delta(j)} \mbox{ for all } j \in \Z/N\Z.
\end{equation}
For any such $\delta$, and any $A \subseteq \Z/N\Z$, the following are
equivalent:
\begin{itemize}
\item[(i)] $\sharp \tilde{g}(A) = \sharp A$;
\item[(ii)] $A$ is a union of cosets of the subgroup $\ell \Z / N
  \Z$ of $\Z/N\Z$;
\item[(iii)] $\sigma \delta(A) = \tilde{g}(A)$.
\end{itemize}
\end{proposition}
\begin{pf}
To prove the first assertion, we exhibit a permutation $\delta$ with
the required property. For $0 \leq i <N$, write
$i=j+c\ell$ with $0 \leq c <m$ and $0 \leq j < \ell$, and set
$\delta(i) = mj + c$. It is routine to verify that $\delta \in S_N$,
and, as
$$ f(I_j) = \bigcup_{d=0}^{m-1} I_{mj+d}, $$
the condition (\ref{delta-prop}) holds.

Now fix a choice of $\delta \in S_N$ satisfying
(\ref{delta-prop}). Since $\sharp \tilde{g}(A) = \sharp \tilde{f}(A)$, the
equivalence of (i) and (ii) follows from Proposition
\ref{size-At}. Since $\sigma \delta \in S_N$, it is immediate that
(iii)$\Rightarrow$(i). It remains to show that (ii)$\Rightarrow$(iii).

Since $f(I_j)=f(I_{j+\ell})$ for each $j$, it follows from
(\ref{delta-prop}) that $\delta$ takes the $m$ elements $j+c \ell$, $0
\leq c <m$ to the $m$ elements $mj+d$, $0 \leq d <m$ in some order.
Thus, if (ii) holds, $\delta$ takes each coset $a +\ell\Z /N\Z$
contained in $A$ to $\tilde{f}(\{a\})$. Thus
$\delta(A)=\tilde{f}(A)$, and applying $\sigma$ gives (iii).
\end{pf}

We consider partitions $\Z / N\Z$ into disjoint non-empty sets: $\Z/ N \Z = A_1
\cup \ldots \cup A_t$. We call the set $\BB=\{A_1, \ldots, A_t\}$ of
subsets of $\Z/N\Z$ a {\em block
  decomposition} of $\Z/N\Z$, and refer to the $A_i$ as {\em blocks}. We say that $\BB$ is {\em trivial} if
$t=1$, and that $\BB$ is $\ell$-stable if, for any $j \in \Z/N\Z$ and
$1 \leq r \leq t$, we have $j \in A_r \Rightarrow j+\ell \in
A_r$. Thus $\BB$ is $\ell$-stable if and only if each $A_r$ is a
union of cosets of the subgroup $\ell \Z / N \Z$ of $\Z/N\Z$. If
$\BB=\{A_1, \ldots, A_t\}$ is a block decomposition and $\sigma \in
S_N$, then $\sigma \BB = \{\sigma(A_1), \ldots, \sigma(A_t)\}$ is also
a block decomposition, and we define the stabiliser $G_\BB$ of $\BB$
as
$$ G_\BB = \{ \sigma \in S_N \ : \  \sigma(\BB)=\BB \} .$$
Then $G_\BB$ is a subgroup of $S_N$.

\begin{lemma} \label{blocks}
Let $f(x)=m x \bmod 1$ and let $N=m \ell$. Let
$\delta$ be as in Proposition \ref{delta-def}. Then, for any $\sigma \in
S_N$, the composite $g= \sigma \circ f$ fails to be mixing if and only
if there is some non-trivial $\ell$-stable block decomposition $\BB$ of
$\Z/N\Z$ such that $\sigma \delta \in G_\BB$.
\end{lemma}
\begin{pf}
Let $\sigma \delta \in G_\BB$ for some non-trivial, $\ell$-stable
block decomposition $\BB$, and let $A$ be a block of $\BB$. Then $A$
is a proper subset of $\Z/N\Z$ which is a union of cosets of
$\ell\Z/N\Z$. Thus $\tilde{g}(A)= \sigma \delta(A)$ by Proposition
\ref{delta-def}, and this set is also a block of $\BB$. Inductively,
we then have $\tilde{g}^r(A)=(\sigma \delta)^r(A)$, and hence
$\sharp \tilde{g}^r(A)=\sharp (\sigma \delta)^r(A)= \sharp A$, for
all $r \geq 0$. It then follows from Lemma \ref{tilde-crit} that $g$
is non-mixing.

Conversely, suppose that $g$ is non-mixing. By Lemma
\ref{tilde-crit}, there is a proper subset $A$ of $\Z/N\Z$ such that
$\sharp \tilde{g}^r(A)=\sharp A$ for all $r\geq 0$. By Proposition
\ref{delta-def} and induction, $\tilde{g}^r(A) =
(\sigma\delta)^r(A)$ for all $r\geq 0$. Moreover, each
$(\sigma\delta)^r(A)$ is a union of cosets of $\ell \Z /N \Z$. Since
$\sigma \delta$ is a permutation, it follows that
$(\sigma\delta)^s(A^c)$ is also a union of cosets for each $s \geq
0$, where $A^c$ is the complement of $A$. Let $\widetilde{\BB}$ be
set of all intersections of the sets $(\sigma\delta)^r(A)$,
$(\sigma\delta)^s(A^c)$ for $r$, $s \geq 0$. Thus $\widetilde{\BB}$
is a collection of subsets of $\Z/N \Z$, each of which is a union of
cosets of $\ell \Z /N \Z$. Let $\BB$ be the collection of minimal
non-empty sets in $\widetilde{\BB}$. Then $\BB$ is an $\ell$-stable
block decomposition and $\sigma \delta \in G_\BB$. Moreover, $\BB$
is non-trivial since $A$ is a union of blocks of $\BB$.
\end{pf}

\begin{remark}
A similar argument shows that $f \circ \sigma$ is non-mixing if and
only if $\delta \sigma \in G_\BB$ for some non-trivial $\ell$-stable
block decomposition.
\end{remark}

\subsection{Asymptotic behaviour as $\ell \to \infty$}
\label{l-asymp}

We continue to assume $N=m\ell$. We shall investigate the proportion
of permutations which do not preserve mixing:
\begin{equation}
p(\ell,m) = \frac{ \sharp\{ \sigma \in S_{m \ell} \ : \ \sigma \circ f
\mbox{ is not
    mixing } \}} { (m\ell)!}.
\end{equation}
By Lemma \ref{blocks}, this is the proportion of permutations such
that $\delta \sigma$ is in the stabiliser of at least one non-trivial
$\ell$-stable block decomposition.

The following Lemma will complete the proof of Theorem \ref{new-mixing-theorem}.

\begin{lemma} \label{asymp}
When $N= m \ell$ with $\ell \geq 6$, we have
$$ p(\ell,m) < 11 \left( \frac{2e}{\ell} \right)^{m-1}. $$
In particular, for each fixed $m \geq 2$ we have $p(\ell,m) \to 0$ as $\ell
\to \infty$.
\end{lemma}

From Lemma \ref{blocks} we have
\begin{equation} \label{sum-B}
  p(\ell,m) \leq {1 \over (m \ell)! } \sum_\BB \sharp G_\BB ,
\end{equation}
where the sum is over all non-trivial $\ell$-stable block decompositions $\BB$.
(This is not an equality since the $G_\BB$ are not disjoint.)
Given integers $1 \leq r_1 \leq \ldots \leq r_j$ with $r_1 + \cdots +
r_j= \ell$,  we consider the contribution to (\ref{sum-B}) from all
block decompositions $\BB$ with block
sizes $mr_1, \ldots, mr_j$. The number of such block decompositions
can be found as follows.  Let us set
$$ n_i(r_1, \ldots, r_j) = \sharp \{ h \ : r_h = i \} $$
and
$$ d(r_1, \ldots, r_j) = \prod_{i=1}^{\ell} n_i(r_1, \ldots, r_j)!. $$
Then the number of $\ell$-stable block decompositions $\BB$ of $\{1, \ldots,
m\ell\}$ with block sizes $mr_1,
\ldots, mr_j$ is
$$ {1 \over d(r_1, \ldots, r_j)} {\ell \choose  r_1, \ldots, r_j
},  $$
where
$$ {\ell \choose r_1, \ldots, r_j } = { \ell! \over r_1 !  \ldots r_j !} $$
is the multinomial coefficient.
Moreover, any such $\BB$ is preserved by a group of permutations
$S_{mr_1} \times \cdots \times S_{m r_j}$ permuting the elements
within each block, but we can also permute the blocks of any given
size amongst themselves. Thus we have
$$ \sharp G_\BB =  d(r_1, \ldots, r_j)
         \left( \prod_{h=1}^j (mr_h)! \right) . $$
The contribution to (\ref{sum-B}) from block decompositions with block
sizes $mr_1, \ldots, mr_j$ is therefore
$$  {1 \over (m \ell)!}
  \left[ d(r_1,\ldots,r_j)  \left( \prod_{h=1}^j (mr_h)! \right)
    \right]
    \left[ {1 \over d(r_1, \ldots, r_j)} {\ell \choose  r_1, \ldots,
      r_j } \right]   $$
which simplifies to
$$   { \ell \choose r_1, \ldots, r_j }
          {m \ell \choose m r_1, \ldots, m r_j }^{-1}.  $$
Thus we may rewrite (\ref{sum-B}) as
\begin{equation} \label{bj-sum}
 p(\ell,m) \leq \sum_{j=2}^{\ell} b_j(\ell),
\end{equation}
where
$$ b_j(\ell) = \sum_{\stackrel{1 \leq r_1 \leq \ldots \leq r_j}{ r_1+ \ldots
  +r_j=\ell}} { \ell \choose  r_1, \ldots, r_j }
         {m \ell \choose m r_1, \ldots, m r_j }^{-1}. $$
The definition of $b_j(\ell)$ makes sense for $j=1$, giving $b_1(\ell)=1$.

\begin{proposition} \label{rec}
For $2 \leq j \leq \ell$, we have
$$ b_j(\ell) \leq  \sum_{r=1}^{\lfloor \ell/j\rfloor} {\ell \choose r} {m \ell
         \choose m r}^{-1} b_{j-1}(\ell - r). $$
\end{proposition}
\begin{pf}
Separating out $r_1$ in the definition of $b_j(\ell)$, we may write
$$  b_j(\ell) \leq \sum_{r_1=1}^{\lfloor \ell/j\rfloor}
    \sum_{ \stackrel{1 \leq r_2 \leq \ldots \leq r_j}{r_2+\cdots + r_j = \ell-r_1}}
       { \ell \choose r_1, \ldots, r_j }
         {m \ell \choose m r_1, \ldots, m r_j }^{-1}. $$
(Note that we have ``$\leq$'' rather than ``='' since the condition $r_2 \geq r_1$
has been weakened to $r_2 \geq
1$.)  The result then follows on using the (easily verified)
identity
$$ {\ell \choose  r_1, \ldots, r_j } = {\ell \choose r_1}
    {\ell -r_1 \choose  r_2, \ldots, r_j }, $$
together with the corresponding identity where all the arguments are
multiplied by $m$.
\end{pf}

\begin{proposition} \label{bj-bound-prop}
Suppose that $m \geq 2$ and $\ell \geq 3$. Then, for $1 \leq j \leq
\ell$, we have
\begin{equation} \label{bj-bound}
 b_j \leq  \left({2e \over \ell}\right)^{(m-1)(j-1)} .
\end{equation}
\end{proposition}
\begin{pf}
We argue by induction on $j$. The result holds for $j=1$ since
$b_1(\ell)=1$.  Suppose that $2\leq j \leq \ell$ and that the result holds
for $j-1$. From Proposition \ref{rec}, we have
\begin{equation} \label{bj-split}
  b_j (\ell) \leq {\ell \choose 1} {m \ell
         \choose m }^{-1} b_{j-1}(\ell - 1) +
\sum_{r=2}^{\lfloor \ell/j\rfloor} {\ell \choose r} {m \ell
         \choose m r}^{-1} b_{j-1}(\ell - r).
\end{equation}
For the first term, we have the estimate
\begin{eqnarray*}
{\ell \choose 1} {m \ell  \choose m}^{-1} b_{j-1}(\ell-1)
  & \leq & \frac{\ell \, (m!)}{(m\ell) (m\ell-1) \ldots
  (m\ell-\ell+1)} \left( \frac{2e}{\ell-1} \right)^{(m-1)(j-2)} \\
 & \leq & \frac{(m-1)!}{m^{m-1} (\ell-1)^{m-1}}
        \left( \frac{2e}{\ell} \cdot
       \frac{\ell}{\ell-1} \right)^{(m-1)(j-2)}  \\
   & = & {1 \over 2} \left( {1 \over 2} \cdot {\ell \over
         \ell-1} \cdot {2 \over \ell} \right)^{m-1}
              \left( \frac{2e}{\ell}
       \cdot \frac{\ell}{\ell-1} \right)^{(m-1)(j-2)}  \\
   & \leq & {1 \over 2^m e^{m-1} }
         \left( \frac{2e}{\ell} \right)^{(m-1)(j-2)}
       \left( \frac{\ell}{\ell-1} \right)^{(m-1)(j-1)}.
\end{eqnarray*}
But
$$ \left( \frac{\ell}{\ell-1} \right)^{(m-1)(j-1)} \leq
  \left( \frac{\ell}{\ell-1} \right)^{(m-1)(\ell-1)} < e^{m-1} $$
since $(1 + {1 \over n})^n$ is an increasing function of $n$ and
$(1 + {1 \over n})^n \to e$ as $n \to \infty$.
As $m \geq 2$, it follows that
\begin{equation} \label{first-term}
{\ell \choose 1}{m \ell \choose m}^{-1} b_{j-1}(\ell-1)
   \leq {1 \over 4 }
         \left( \frac{2e}{\ell} \right)^{(m-1)(j-1)}.
\end{equation}

We now consider each term in the sum in \eqref{bj-split}.
For $2 \leq r \leq \lfloor \ell/j \rfloor$, we have
$$ {\ell  \choose  r}^m \leq {m \ell  \choose m r}. $$
(This is obvious combinatorially: some of the ways of choosing $m r$
objects from $m \ell$ are given by choosing $r$ objects from the first
$\ell$, then another $r$ from the second $\ell$, and so on.) Also,
since $2 \leq r \leq \ell/2$, we have
$$ {\ell \choose 2 } \leq {\ell \choose r}. $$
Thus
$$ {\ell \choose r} {m \ell\choose m r}^{-1} \leq {\ell \choose
  r}^{1-m} \leq {\ell \choose 2}^{1-m} =
   {2^{m-1} \over \ell^{m-1} (\ell-1)^{m-1} } . $$
From the induction hypothesis, we have
$$ b_{j-1}(\ell-r) \leq \left(\frac{2e}{\ell-r}\right)^{(m-1)(j-2)} \leq
   \left(\frac{2e}{\ell}\right)^{(m-1)(j-2)}
     \left(\frac{j}{j-1}\right)^{(m-1)(j-2)},  $$
and $\big(j/(j-1) \big)^{j-2}<\big(j/(j-1) \big)^{j-1}<e$.
Thus
\begin{eqnarray*}
\sum_{r=2}^{\lfloor \ell/j\rfloor} {\ell \choose r} {m \ell
         \choose m r}^{-1} b_{j-1}(\ell - r)
 & < & {\ell \over j} {2^{m-1} \over \ell^{m-1} (\ell-1)^{m-1} }
            \left(\frac{2e}{\ell}\right)^{(m-1)(j-2)} e^{m-1} \\
 & = & {\ell \over j(\ell-1)^{m-1} }
           \left(\frac{2e}{\ell}\right)^{(m-1)(j-1)}.
\end{eqnarray*}
But as $j \geq 2$, $m\geq 2$ and $\ell \geq 3$, we have
$$ {\ell \over j(\ell-1)^{m-1} }  \leq {\ell \over 2(\ell-1)} \leq
{3 \over 4}. $$
Substituting the last estimate and \eqref{first-term} into
\eqref{bj-split}, we therefore obtain
$$ b_j \leq  \left({2e \over \ell}\right)^{(m-1)(j-1)} , $$
which completes the induction.
\end{pf}

\begin{pf-asymp}
Since we assuming $\ell \geq 6$, we have $2e/\ell<1$. It then follows
from \eqref{bj-sum} and Proposition \ref{bj-bound-prop} that
\begin{eqnarray*}
  p(\ell,m) & < & \sum_{j=2}^{\infty} \left( {2e \over \ell}
\right)^{(m-1)(j-1)}  \\
  & = &    \left( {2e \over \ell} \right)^{m-1}
  \left[1-\left( {2e \over \ell} \right)^{m-1} \right]^{-1}.
\end{eqnarray*}
As $m \geq 2$ and $2e/\ell < \frac{10}{11}$, this gives
$$ p(\ell,m)   <
       \left[1-\left( {2e \over \ell} \right) \right]^{-1}
        \left( {2e \over \ell} \right)^{m-1} <
    11 \left( {2e \over \ell} \right)^{(m-1)}, $$
as required.
\end{pf-asymp}

\subsection{The proportion of non-mixing permutations}
\label{l-small}

In this section, we will use
the Inclusion-Exclusion Principle (see e.g.~\cite[p.~21]{Slomson})
to give explicit formulae for the proportion $p(\ell,m)$ of non-mixing
permutations when $N=m \ell$ with $\ell$ small.


The stabiliser of any non-trivial $\ell$-stable block decomposition
contains the subgroup $H \cong S_m \times \ldots \times S_m$ of order
$(m!)^\ell$ which permutes the $m$ elements of each coset amongst
themselves.
In order to refer to specific block decompositions, we let
$C_1, \ldots, C_\ell$ denote the cosets of $\ell \Z /N \Z$
in $\Z/N\Z$ (in some order). Giving an $\ell$-stable block
decomposition amounts to giving a partition of $\{C_1, \ldots,
C_\ell\}$, and we denote the block decomposition by the corresponding
partition of the set of indices $\{1, \ldots, \ell\}$.
Thus $\{1, \ldots, \ell-1\}, \{\ell\}$ represents
the $\ell$-stable block decomposition consisting of the two blocks
$C_1 \cup \ldots \cup C_{\ell-1}$ of size $(\ell-1) m$ and $C_\ell$ of
size $m$.

%
%
\subsubsection{$\ell=2$}  \label{l2}

There is only one non-trivial $\ell$-stable block
decomposition. This has two blocks, each of size $m$. Its stabiliser
contains $H$ and also contains elements swapping the two blocks, so has order $2\sharp H$. Thus
$$ p(2,m) = { 2\sharp H \over (2m)!} = {2m-1 \choose m}^{-1}. $$
In particular, taking $m=2$, we get $p(2,2) = 1/3$. Thus, when the
doubling map $f(x)=2x \bmod 1$ is composed with permutations $\sigma$ of the 4 equal subintervals of $[0,1)$,
those $\sigma \in S_4$ for which $f \circ \sigma$ is not
mixing form a single coset of a subgroup of index 3 in $S_4$. (Any such
subgroup is dihedral of order 8.)

\subsubsection{$\ell = 3 $} \label{l3}

There are 4 non-trivial $\ell$-stable block decompositions:
$$ \mathrm{(i)} \;  \{ 1, 2\}, \; \{3\}; \quad
  \mathrm{(ii)} \;  \{ 1, 3\}, \; \{2\}; \quad
  \mathrm{(iii)} \;  \{ 2, 3\}, \; \{1\}; \quad
  \mathrm{(iv}) \;  \{ 1 \}, \; \{ 2 \},\; \{3\}.  $$
%
The stabiliser of any one of the block decompositions (i), (ii), (iii)
has order $(2m)! m! = {2m \choose m} \sharp H$ since it contains any
permutation of the $2m$ elements in the block consisting of 2 cosets.
The stabiliser of the block decomposition (iv) has order $6 \sharp H$ since
we may permute the 3 blocks amongst themselves in $3!=6$ ways.

We now consider the stabilisers of any of the ${4 \choose 2}=6$ pairs
of the block decompositions. First consider the 3 pairs
consisting of any two of (i), (ii) or (iii). Any permutation fixing such a
pair must fix each coset, so the stabiliser of any of these 3
pairs is just $H$. A permutation stabilising (say) (i) and (iv)
could also swap the cosets $C_1$ and $C_2$, so the stabilisers of the
other 3 pairs have orders $2 \sharp H$. The stabiliser of any 3 (or all
4) block decompositions is again just $H$. Thus the precise number of
permutations in $S_{3m}$ fixing at least one of the block
decompositions is
$$ \left( 3 {2m \choose m} + 6 - 3 - 3\times 2 + {4 \choose 3} -1
\right) \sharp H = 3 {2m \choose m} \sharp H  = 3 (2m)! m!.$$
Hence
$$ p(3,m) = {1 \over (3m)!} \times 3 (2m)! m! = {3m-1 \choose 2m}^{-1}. $$
In particular, $p(3,2) = 1/5$.

\subsubsection{$\ell=4$} \label{l4}

There are 14 non-trivial $\ell$-stable block
decompositions, but we only need to consider the
4 block decompositions with block sizes 3, 1 and the 3 block
decompositions with block sizes 2, 2, since any permutation
stabilising a non-trivial block decomposition must stabilise one of
these. We can then apply to the Inclusion-Exclusion Principle to the
stabilisers of these 7 block decompositions, by considering all
possible pairs, and, for each pair, considering any ways of extending
the pair to a larger subset of the blocks with
stabiliser larger than $H$. After some simplification, we
obtain the formula
$$ p(4,m)  = \left[
  4{3m \choose m, m, m} + 6{2m \choose m}^2  -12{2m \choose m} \right]
  \frac{ (m!)^4}{(4m)!}  . $$
In particular, we find
$$  p(4,2) = {1 \over 5}, \qquad
     p(4,3) = {37 \over 1540}. $$
Note that, in contrast to the cases $\ell=2$ and $\ell=3$,
$p(4,m)$ is not in general the reciprocal of an integer.

\section{Mixing rates for $mx \bmod 1$}\label{sec_proof}

In this section, we prove Theorem \ref{thm:perm-mix}.
The computation of the essential spectral radius $r_{ess}$ for $\sigma
\circ f$ is straightforward since $\sigma\circ f$ is piecewise linear with
constant slope $1/m$. Hence Theorem \ref{thm:perm-mix}(i)
is a consequence of \cite{Keller84}.


We now turn to Theorem \ref{thm:perm-mix}(ii). This requires a detailed study of
the eigenvalues of the Fredholm matrices $\Phi(z)$ associated to
$\sigma\circ f$. We first give the required background on Fredhom
matrices, see \cite{Mori89,Mori98}.

\subsection{Fredholm matrices}

Consider a piecewise linear Markov map $f:I\to I$, with finite
partition $\mathcal{P}=\{I_i\}_{i=1}^{q}$, and representative
transition matrix $B$. Here $B$ is a $q\times q$ matrix with
$B_{ij}=1$ if $I_j\subset f(I_i)$, and $B_{ij}=0$ if $f(I_i)\cap
I_j=\emptyset$. We will assume that $f$ is differentiable on the
interior of each element of $\mathcal{P}$. If $\mathcal{L}_f$ is the
transfer operator, and $J\subset I$, we consider the power series
defined on $\mathbb{C}\times D$, with $D\subset\mathbb{C}$:
\begin{equation}
s^{J}(z,x)=\sum_{n=0}^{\infty}z^n\mathcal{L}^{n}_{f}(\mathcal{X}_{J})(x)=\mathcal{X}_{J}(x)+\sum_{n=1}^{\infty}z^n\mathcal{L}^{n}_f
(\mathcal{X}_{J})(x),
\end{equation}
where $\mathcal{X}_J(x)$ is the indicator function of $J$. When
$J=I_i\in\mathcal{P}$, we will write
$s^J(z,x)$ as $s^{(i)}(z,x).$ We let $\underline{s}(z,x)$ be the vector $(s^{(i)}(z,x))_{i=1}^{q}$, and
similarly $\underline{\mathcal{X}}(x)=(\mathcal{X}_{(i)}(x))_{i=1}^{q}$.
For a Markov system we have the following result.
\begin{lemma}\label{prop:fredholm}
For a piecewise linear Markov map $f:I\to I$ with finite partition $\mathcal{P}=\{I_i\}_{i=1}^{q}$,
there exists a $q\times q$ matrix $\Phi(z)$, and such that
\begin{equation}
\underline{s}(z,x)=(I-\Phi(z))^{-1}\underline{\mathcal{X}}(x).
\end{equation}
\end{lemma}
The matrix $\Phi(z)$ in Lemma \ref{prop:fredholm} is called a
\emph{Fredholm matrix}.

\noindent \begin{pf}
We will consider the Markov case where the slope is constant on each
$I_i$ (but not constant globally). Our proof is a slight adaption of
the calculations in \cite{Mori89, Mori98}. In particular we will
obtain an explicit form of $\Phi(z)$. First of all, by definition of
$\mathcal{L}_f$ we have
\begin{equation*}
s^{J}(z,x)=\mathcal{X}_{J}(x)+\sum_{n=1}^{\infty}z^n\sum_{f^n(y)=x}\frac{\mathcal{X}_{J}(y)}{|(f^n)'(y)|}.
\end{equation*}
If $J=I_i\in\mathcal{P}$, the following hold:
\begin{equation*}
\begin{split}
s^{(i)}(z,x) &=\mathcal{X}_{(i)}(x)+z\sum_{n=1}^{\infty}z^{n-1}\sum_{f^{n-1}(f(y))=x}
\frac{\mathcal{X}_{(i)}(y)}{|(f^{n-1})'(f(y))f'(y)|}\\
&=\mathcal{X}_{(i)}(x)+z\sum_{n=1}^{\infty}\frac{z^{n-1}}{|(f'\mid I_i)|}\sum_{f^{n-1}(\tilde{y})=x}\sum_{\underset{B_{ij}=1}{i,j}}
\frac{\mathcal{X}_{(j)}(\tilde{y})}{|(f^{n-1})'(\tilde{y})|}\\
&=\mathcal{X}_{(i)}(x)+\left(\frac{z}{|(f'\mid I_i)|}\sum_jB_{ij}\right)s^{(j)}(z,x).
\end{split}
\end{equation*}
Hence we obtain a $q\times q$ matrix $\Phi(z)$, with $\Phi(z)_{ij}=\{z/|(f'\mid I_i)|\}B_{ij},$ and
\begin{equation}
\underline{s}(z,x)=(I-\Phi(z))^{-1}\underline{\mathcal{X}}(x).
\end{equation}
This completes the proof.
\end{pf}

Given the Fredholm matrix $\Phi(z)$, we define the \emph{Fredholm
determinant} to be the quantity $D(z)=\det(I-\Phi(z)).$ For
piecewise-linear expanding (Markov) systems, the Fredholm matrix and
Fredholm determinant have the following properties (see
\cite{Mori89, Mori98}) which are useful in the sequel:
\begin{enumerate}
\item The number of ergodic components of $f$ is equal to the dimension of the eigenspace of $I-\Phi(1)$ associated to the
eigenvalue of value zero. The number of ergodic components is also equal to the order of the zero at $z=1$ in the equation
$\det(I-\Phi(z))=0$.
\item If zero is a simple eigenvalue of $I-\Phi(1)$ then the system is ergodic. Moreover if
$\{|z|=1\}\cap\mathrm{Spec}(\mathcal{L}_{f}|_{\operatorname{BV}})=\{1\}$
then the system is mixing
\item If $\lambda\in\mathbb{C}$ and $|\lambda|>r_{ess}$, then $\lambda\in\text{Spec}
(\mathcal{L}_{f}|_{\operatorname{BV}})$ if and only if
$z=\lambda^{-1}$ is a zero of $D(z)$, i.e. $D(1/\lambda)=0$.
\item If $D(1/\lambda)=0$ then $\lambda$ is an eigenvalue
of $\mathcal{L}_{f}|_{\operatorname{BV}}$.
\end{enumerate}

\subsection{Computation of Fredholm matrix eigenvalues}

We now consider Fredholm matrices for our maps $\sigma \circ f$ with
$f(x)=mx \bmod 1$ and $\sigma \in S_N$, where we assume that $N>m$ and
$\gcd(m,N)=1$. These matrices are attached
to a partition of $[0,1]$ on which $\sigma \circ f$ is Markov, so we
first need to determine such a partition.
For $k \geq 1$, consider the partition
$$\mathcal{P}_k:=\left\{ [(j/k,(j+1)/k) \ : \ 0 \leq j \leq k-1   \right\}$$
of $[0,1)$ into $k$ equal subintervals. Then the map $f$ is Markov
w.r.t.~$\mathcal{P}_{m}$, while the map $\sigma$ is Markov
w.r.t.~$\mathcal{P}_N$. The map $\sigma \circ f$, however, is in
general not Markov w.r.t.~either of these partitions. For example
consider $m=2, N=3$. Clearly any $\sigma \in S_3$ is Markov on the
partition
$$\mathcal{P}_3=\{[0,1/3),[1/3,2/3),[2/3,1)\}.$$
However, if we take the permutation $\sigma$ interchanging the last
two subintervals, then we have
$$ \sigma \circ f(x) = \begin{cases}
     2x & \mbox{if } 0 \leq x < 1/6,  \cr
     2x + 1/3 & \mbox{if } 1/6 \leq x <1/3, \end{cases} $$
so that $\sigma \circ f$ is not continuous on $[0,1/3]$ and hence not
Markov on $\mathcal{P}$. In general, to ensure that $\sigma \circ f$
is Markov for all $\sigma \in S_N$, we must work with the partition
$\mathcal{P}_{Nm}$.

Due to our specific choice $f(x) = mx \bmod 1$, the $Nm \times Nm$
matrix $\Phi(1)$ is precisely the probability-transition matrix
between the Markov states, and has all its entries in $\{0, 1/m\}$.
If $\lambda\in\mathrm{Spec}(\mathcal{L}_{f}|_{\operatorname{BV}})$
then we know that $z=1/\lambda$ is a solution to
$D(z)=\det(I-\Phi(z))=0$. It is therefore an equivalent problem to
consider the corresponding equation (in $\lambda$) to
$\det(B-\lambda I)$=0, where $B$ is the state transition matrix
(with entries in $\{0,1\}$). Hence if $\tilde\lambda$ is an
eigenvalue of $B$, then
$\lambda=\tilde\lambda/m\in\mathrm{Spec}(\mathcal{L}_{f}|_{\operatorname{BV}}).$

Note that in our case, $m \Phi(1)$ is precisely the state transition matrix $B$.
We will show that the eigenvalues of $\Phi(1)$
can in fact be determined from the $N \times N$ transition matrix
associated with the partition $\mathcal{P}_N$ (on which $\sigma \circ f$
need not be Markov).

We must first define some notation. Following the conventions of
Section \ref{sec_mixvsnonmix}, we index the subintervals in
$\mathcal{P}_k$ by $\{0, 1, \ldots, k-1\}$. We therefore begin the
numbering of the rows and columns in the associated matrices from
$0$. We define $A(m,N)$ and $B(m,N)$ to be the state transition
matrices for $f$ w.r.t.~$\mathcal{P}_N$ and $\mathcal{P}_{Nm}$
respectively. Thus for $0 \leq i, j \leq N-1$ we have
$$ A(m,N)_{ij} = \begin{cases} 1 & \mbox{if } j \equiv mi+d \bmod N
  \mbox{ with } 0
  \leq d \leq m-1, \cr
   0 & \mbox{otherwise,} \end{cases} $$
and for $0 \leq i, j \leq Nm-1$ we have
$$ B(m,N)_{ij} = \begin{cases} 1 & \mbox{if } j \equiv mi+d \bmod Nm
                      \mbox{ with } 0 \leq d \leq m-1, \cr
   0 & \mbox{otherwise}. \end{cases}$$
For example, when $m=2$ and $N=3$, we have
\begin{equation*}
A(2,3)=\begin{pmatrix}
1 & 1 & 0\\
1 & 0 & 1 \\
0 & 1 & 1
\end{pmatrix},
\quad
B(2,3)=\left( \begin{array}{cc|cc|cc}
1 & 1 & 0 & 0 & 0 & 0\\
0 & 0 & 1 & 1 & 0 & 0\\
0& 0 & 0 & 0 & 1 & 1 \\ \hline
1 & 1 & 0 & 0 & 0 & 0\\
0 & 0 & 1 & 1 & 0 & 0\\
0& 0 & 0 & 0 & 1 & 1 \\
\end{array}\right).
\end{equation*}
The eigenvalues for $A(2,3)$ are $\{\pm1,2\}$, while those for
$B(2,3)$ are $\{\pm 1,2,0\}$, where the eigenspace for the eigenvalue
$0$ has dimension $3$. In the case
$m=3,N=5$ we have:
\begin{equation*}
A(3,5)=\begin{pmatrix}
1 & 1 & 1 & 0 & 0 \\
1 & 0 & 0 & 1 & 1 \\
0 & 1 & 1 & 1 & 0 \\
1 & 1 & 0 & 0 & 1\\
0 & 0 & 1 & 1 & 1\\
\end{pmatrix},
\end{equation*}
and the eigenvalues for $A(3,5)$ are $\{3,\pm i,\pm 1\}$.
Note that all row sums and columns sums in both $A(m,N)$ and $B(m,N)$
are $m$. Each row in either matrix consists of $m$ consecutive
occurrences of $1$ (where, in the case of $A(m,N)$ these may ``wrap
around'' from the last column to the first). The rows of $B(m,N)$
naturally fall into $m$ identical blocks each consisting of $N$ rows,
and the columns into $N$ blocks each consisting of $m$ identical
columns, as indicated for $B(2,3)$ above.

The corresponding state transition matrices for $\sigma \circ f$ are
obtained by permuting the columns of $A(m,N)$ and $B(m,N)$. More
precisely, given a permutation $\sigma$ of $\{0, \ldots, N-1\}$, let
$P(\sigma)$ be the $N \times N$ permutation matrix given by
$$ P(\sigma)_{ij} = \begin{cases} 1 & \mbox{ if } j = \sigma(i), \cr
                                 0 & \mbox{otherwise,} \end{cases} $$
and let $Q(\sigma)$ be the $Nm \times Nm$ matrix obtained by replacing
each entry $1$ (respectively, $0$) in $P(\sigma)$ by an $m \times m$
identity matrix (respectively, zero matrix). Then the state transition
matrices for $\sigma \circ f$ w.r.t.~the partitions $\mathcal{P}_N$
and $\mathcal{P}_{Nm}$ are $A(m,N)P(\sigma)$ and $B(m,N)Q(\sigma)$
respectively. For example, if $m=2$, $N=3$ and $\sigma$ is the 3-cycle
$(0,1,2)$ then
\begin{equation*}
P(\sigma)=\begin{pmatrix}
0 & 1 & 0\\
0 & 0 & 1 \\
1 & 0 & 0
\end{pmatrix},
\quad
Q(\sigma)=\left( \begin{array}{cc|cc|cc}
0 & 0 & 1 & 0 & 0 & 0\\
0 & 0 & 0 & 1 & 0 & 0\\ \hline
0 & 0 & 0 & 0 & 1 & 0 \\
0 & 0 & 0 & 0 & 0 & 1\\ \hline
1 & 0 & 0 & 0 & 0 & 0\\
0 & 1 & 0 & 0 & 0 & 0 \\
\end{array}\right),
\end{equation*}
so that
\begin{equation*}
A(2,3)P(\sigma)=\begin{pmatrix}
0 & 1 & 1\\
1 & 1 & 0 \\
1 & 0 & 1
\end{pmatrix},
\quad
B(2,3)Q(\sigma)=\left( \begin{array}{cc|cc|cc}
0 & 0 & 1 & 1 & 0 & 0\\
0 & 0 & 0 & 0 & 1 & 1\\ \hline
1 & 1 & 0 & 0 & 0 & 0 \\
0 & 0 & 1 & 1 & 0 & 0\\ \hline
0 & 0 & 0 & 0 & 1 & 1\\
1 & 1 & 0 & 0 & 0 & 0 \\
\end{array}\right).
\end{equation*}
Note also that
$P(\sigma) A(m,N)$ is the matrix obtained by applying the {\em
inverse} permutation $\sigma^{-1}$ to the {\em rows} of $A(m,N)$.

To determine the mixing rate of $\sigma \circ f$, we need to
investigate the eigenvalues of the Fredholm matrix $\Phi(1)=m^{-1}
B(m,N)Q(\sigma)$ corresponding to the partition $\mathcal{P}_{Nm}$ on which
$\sigma \circ f$ is Markov. Clearly $\lambda$ is an eigenvalue of
$\Phi(1)$ if and only if $m \lambda$ is an eigenvalue of $B(m,N) Q(\sigma)$, so
it suffices to find the eigenvalues of the latter $Nm \times Nm$
matrix. In fact we only need consider $N \times N$ matrices.

\begin{lemma} \label{shrink-matrix}
For all $m$, $N$ and all $\sigma \in S_N$, the nonzero eigenvalues of
$B(m,N)Q(\sigma)$ are the same as those of $A(m,N)P(\sigma) $.
\end{lemma}
\begin{pf}
For brevity, we write $A=A(m,N)$, $B=B(m,N)$, $P=P(\sigma)$ and
$Q=Q(\sigma)$.

We view $BQ$ as determining a linear endomorphism $\theta$ on the space
$V=\C^{Nm}$ of column vectors. Clearly $BQ$ has rank $N$, since the
first $N$ rows are linearly independent and the remaining rows merely
repeat these. The kernel $W$ of $\theta$ therefore has dimension
$N(m-1)$, and $\theta$ induces an endomorphism $\overline{\theta}$ on
the quotient space $V/W$ of dimension $N$. The eigenvalues of
$\theta$ (that is, of $BQ$) are therefore the eigenvalues of
$\overline{\theta}$, together with the eigenvalue $0$ of multiplicity
$N(m-1)$ coming from $W$. The result will therefore follow if we show
that the matrix $AP$ represents $\overline{\theta}$.

We define vectors $\mathbf{v}^{r,s}$ for $0 \leq r \leq N-1$, $0 \leq s\leq m-1$
(independent of $\sigma$) as follows. For $s=0$, set
$$ \mathbf{v}^{r,0}_i = \begin{cases}
  1 & \mbox{if } i=mr, \cr
  0 & \mbox{otherwise}, \end{cases} $$
and for $s>0$,
$$ \mathbf{v}^{r,s}_i = \begin{cases}
  -1 & \mbox{if } i=mr, \cr
  1 & \mbox{if } i=mr+s, \cr
  0 & \mbox{otherwise}. \end{cases} $$
For example, if $m=2$ and $N=3$ we have
{\scriptsize
$$
 \mathbf{v}^{0,0} = \left(\begin{array}{r}
      1 \\ 0 \\ \hline 0 \\ 0 \\ \hline 0 \\ 0 \end{array}\right),
\quad
    \mathbf{v}^{0,1} = \left(\begin{array}{r}
      -1 \\ 1 \\ \hline 0 \\ 0 \\ \hline 0 \\ 0 \end{array}\right),
\quad
\mathbf{v}^{1,0} = \left(\begin{array}{r}
      0 \\ 0 \\ \hline 1 \\ 0 \\ \hline 0 \\ 0 \end{array}\right),
\quad
    \mathbf{v}^{1,1} = \left(\begin{array}{r}
      0  \\ 0 \\ \hline -1 \\ 1 \\ \hline 0 \\ 0 \end{array}\right),
\quad
\mathbf{v}^{2,0} = \left(\begin{array}{r}
      0 \\ 0 \\ \hline 0 \\ 0 \\ \hline 1 \\ 0 \end{array}\right),
\quad
    \mathbf{v}^{2,1} = \left(\begin{array}{r}
      0 \\ 0 \\ \hline 0 \\ 0 \\ \hline 1 \\ -1 \end{array}\right),
$$}
where the horizontal lines correspond to the division of the columns of
$B(2,3)Q(\sigma)$ into blocks.

It is clear that the $\mathbf{v}^{r,s}$ form a basis for $V$,
and that if $s \neq 0$ then $BQ\mathbf{v}^{r,s}=0$. Hence the $N(m-1)$
vectors $\mathbf{v}^{r,s}$ for $s \neq 0$ form a basis for $W$.
Thus the $N$ cosets $\mathbf{v}^{r,0} + W$ form a basis for $V/W$. If
we partition $BQ$ into $m \times m$ blocks (as in the above example), the
matrix of $\overline{\theta}$ with respect to this basis is then obtained
by replacing each block with the sum of one of its (identical)
columns. This gives precisely the matrix $AP$.
\end{pf}

We next consider a matrix related to $A(m,N)$ but with eigenvalues
that are easy to determine. By permuting the {\em rows} of $A(m,N)$,
we can obtain a symmetric circulant matrix $C(m,N)$. Its explicit
description depends on the parity of $m$. Let
$$ \delta = \begin{cases} (1-m)/2 & \mbox{if } m \mbox{ is odd;} \cr
                          (1-m+N)/2 & \mbox{if } m \mbox{ is even.}
            \end{cases} $$
Then  $\delta \in \Z$ in both cases since $\gcd(m,N)=1$, and
$C(m,N)$ has entries
\begin{equation}\label{equ_circulant}
 C(m,N)_{ij} = \begin{cases} 1 & \mbox{if } j \equiv i+ \delta +r \bmod N
  \mbox{ with } 0 \leq r \leq m-1, \cr
    0 & \mbox{otherwise} \end{cases}
\end{equation}
for $0 \leq i, j  \leq N-1$. Observe that $C(m,N)$ is indeed symmetric since
$$ j \equiv i+ \delta +r \bmod N \LRA i \equiv j + \delta + (m-1-r) \bmod N.$$
For example,
$$ C(2,5) = \begin{pmatrix} 0 & 0 & 1 & 1 & 0 \\
                            0 & 0 & 0 & 1 & 1 \\
                            1 & 0 & 0 & 0 & 1 \\
                            1 & 1 & 0 & 0 & 0 \\
                            0 & 1 & 1 & 0 & 0 \end{pmatrix}; \qquad
   C(3,5) = \begin{pmatrix} 1 & 1 & 0 & 0 & 1 \\
                            1 & 1 & 1 & 0 & 0 \\
                            0 & 1 & 1 & 1 & 0 \\
                            0 & 0 & 1 & 1 & 1 \\
                            1 & 0 & 0 & 1 & 1 \end{pmatrix}. $$

Since $C(m,N)$ is a real symmetric matrix, its eigenvalues are
real. Since $C(m,N)$ is a circulant matrix, we can write these
eigenvalues down explicitly. Let $\omega_j=e^{2\pi ij/N}$ for $0 \leq
j <N$, and let
\begin{equation} \label{def-vj}
 \mathbf{v}_j=\left(1, \omega_j,
 \omega_j^2,\ldots,\omega_j^{N-1}\right)^{T}.
\end{equation}
Then $\mathbf{v}_j$ is an eigenvector for $C(m,N)$ with eigenvalue
$$ \lambda_j = \sum_{r=0}^{m-1} \omega_j^{\delta + r}. $$
Although the $\lambda_j$ are not necessarily distinct, the $N$
eigenvectors $\mathbf{v}_j$ are linearly independent since
$$ \det( \omega_j^k)_{0 \leq j,k <N} = \prod_{j<k} (\omega_k-\omega_j)
\neq 0, $$
so there are no further eigenvalues. Trivially
$\lambda_0=m$. For $j \neq 0$, we have
$$ \lambda_j = \frac{\omega_j^\delta (\omega_j^m-1)}{\omega_j-1}. $$
Writing $\zeta_j=e^{\pi i j/N}$, so that $\omega_j=\zeta_j^2$, we then
  have
\begin{equation} \label{lambda-val}
 \lambda_j = \zeta_j^{2\delta +m-1}
   \left( \frac{\zeta_j^m-\zeta_j^{-m}}{\zeta_j-\zeta_j^{-1}} \right)
   = (-1)^{(m-1)j} \frac{ \sin( m j \pi /N)}{\sin(j\pi / N)},
\end{equation}
since $2\delta +m-1= 0$ (resp.\ $N$) if $m$ is odd (resp.\ even).
In particular,
\begin{equation} \label{C-det}
 \det(C)  = \prod_{j=0}^{N-1}  \lambda_j
  = \pm  m \prod_{j=1}^{N-1}
    \frac{\sin(mj\pi/N)}{\sin(j\pi/N)} = \pm m,
\end{equation}
using the fact that the residues $mj \bmod N$ are just the residues
$j \bmod N$ in some order because $\gcd(m,N)=1$.
It follows easily from (\ref{lambda-val}) that, for $j \neq 0$, we have
$$ \lambda_{N-j} = (-1)^{(N-1)(m-1)} \lambda_j = \lambda_j, $$
where the second equality holds since $N$ and $m$ cannot both be even.
%
%
%

We also mention some variants of $C(m,N)$. Firstly,
cyclically permuting the rows of $C(m,N)$ gives circulant matrices
$C^{(h)}(m,N)$ for which the $\mathbf{v}_j$ are eigenvectors with
eigenvalues $\omega_j^h \lambda_j$. Secondly, we may permute the rows
of $C(m,N)$ to obtain the anticirculant matrix $C'(m,N)$ with the same
first row as $C(m,N)$; for example
$$ C'(2,5) = \begin{pmatrix} 0 & 0 & 1 & 1 & 0 \\
                            0 & 1 & 1 & 0 & 0 \\
                            1 & 1 & 0 & 0 & 0 \\
                            1 & 0 & 0 & 0 & 1 \\
                            0 & 0 & 0 & 1 & 1 \end{pmatrix}; \qquad
   C'(3,5) = \begin{pmatrix} 1 & 1 & 0 & 0 & 1 \\
                            1 & 0 & 0 & 1 & 1 \\
                            0 & 0 & 1 & 1 & 1 \\
                            0 & 1 & 1 & 1 & 0 \\
                            1 & 1 & 1 & 0 & 0 \end{pmatrix}. $$
Explicitly, the entries of $C'(m,N)$ are
$$  C'(m,N)_{ij} = \begin{cases} 1 & \mbox{if } j \equiv -i+ \delta +r \bmod N
  \mbox{ with } 0 \leq r \leq m-1, \cr
    0 & \mbox{otherwise} \end{cases} $$
for $0 \leq i, j  \leq N-1$. The eigenvalues of $C'(m,N)$ are real
since any anticirculant matrix is symmetric. Using \cite[Theorem
  2]{Chao}, we can write down these eigenvalues explicitly: they are
$m$, $\lambda_{N/2}$ (for $N$ even), and both values of $\pm
 \sqrt{\lambda_j \lambda_{N-j}}=\pm \lambda_j$ for $j \neq 0$, $N/2$.

The maximum value of $|\lambda_j|$ for $1 \leq j \leq N-1$ is attained
at $j=1$. Although this is essentially elementary, it is trickier to
verify than it might appear, so we include a proof.

\begin{proposition}\label{max-circ-eval}
$$ \max_{1 \leq j \leq N-1} \left|\frac{\sin(mj\pi/N)}{\sin(j\pi/N)}\right| =\frac{\sin(m\pi/N)}{\sin(\pi/N)}. $$
\end{proposition}
\begin{pf}
Since $|\sin(\pi k \pm x)|=|\sin x|$ for all $k \in \Z$, we may assume
that $1 \leq m \leq N/2$, and moreover it suffices to take $1 \leq j
\leq N/2$. We consider the two functions $u(x)=\sin mx / \sin x $ and
$v(x)=1/\sin x$ on the interval $(0,\pi)$. Now $u(x)$ has precisely
$m-1$ zeros on this interval, at $x= h \pi /m$ for $1 \leq h \leq
m-1$. Since $u(x)$ may be written as a polynomial of degree $m-1$
in $\cos x$, and $\cos x$ is monotonically decreasing on this
interval, it follows that $u(x)$ has precisely $m-2$ stationary
points, one in each of the intervals $(h\pi/ m, (h+1)\pi/m)$ for $1
\leq h \leq m-2$. In particular, as $\lim_{x \to 0} u(x)=m$, it
follows that $u(x)$ is positive and decreasing on $(0,\pi/m)$, so that
$u(\pi/N)>u(j\pi/N)\geq 0$ if $2 \leq j \leq N/m$.  On the other hand,
as $v(x)$ is positive and decreasing throughout $(0,\pi/2)$, we have
for $N/m \leq j\leq N/2$ that $|u(j \pi/N)| \leq v(j \pi/N) \leq
v(\pi/m) < v(\pi/2m)$.  But $v(\pi/2m)=u(\pi/2m)\leq u(\pi/N)$ as $m
\leq N/2$. Hence $|u(j\pi/N)| <u(\pi/N)$ for $2 \leq j \leq N/2$, as
required.
\end{pf}

We now seek to relate the eigenvalues of the matrices
$A(m,N)P(\sigma)$ to those of $C(m,N)$. After scaling by $1/m$,
these matrices become doubly stochastic.
Our next result gives some information on the behaviour of the
eigenvalues of a column stochastic matrix under permutation of its
columns (or, more generally, under right multiplication by an
orthogonal, column stochastic matrix).

Recall that an $N \times N$ matrix is {\em row} (respectively, {\em
  column}) {\em stochastic} if its entries are non-negative real
numbers and the sum of each row (respectively, column) is $1$.  It is
{\em doubly stochastic} if it is both row and column stochastic.  The
product of two row (respectively, column, doubly) stochastic matrices
is again row (respectively, column, doubly) stochastic. For $A(m,N)$
as above, the probability transition matrices $m^{-1} A(m,N)$ are
doubly stochastic. Any permutation matrix $P(\sigma)$ is doubly
stochastic and orthogonal.

We view our matrices as linear maps on the
space $\C^N$ of column vectors, endowed with the usual complex inner
product $(\xx,\yy)=\sum_{j=1}^N x_j \overline{y}_j$ for $\xx=(x_1,
\ldots, x_N)^T$, $\yy=(y_1, \ldots, y_N)^T$, and we write
$||\xx||=\sqrt{(\xx,\xx)}$ for $\xx \in \C^N$.
Any row stochastic matrix has the obvious eigenvector $\ee=(1,\ldots,1)^T$
with eigenvalue 1. It is well-known that any eigenvalue $\lambda$
satisfies $|\lambda| \leq 1$. If $B$ is a column stochastic matrix
then $\ee$ is not necessarily an eigenvector for $B$, but if
$(\xx,\ee)=0$ then $(B\xx,\ee)=0$, so that $B$ preserves the subspace
$V_0$ of vectors in $\C^N$ perpendicular to $\ee$.

\begin{lemma} \label{circulant-worst}
Let $B$ be an $N \times N$ column stochastic matrix. Then the
eigenvalues of $B^T B$ on $V_0$ are real and non-negative. Let $\eta$
be the largest of these, and let $P$ be an $N \times N$ orthogonal,
column stochastic matrix (e.g.~a permutation matrix).  Then every
eigenvalue $\lambda$ of $BP$ on $V_0$ satisfies
$$ |\lambda| \leq \sqrt{\eta}. $$
Moreover, if $B$ is a circulant matrix then
$$ \sqrt{\eta} = \max\{ \ |\lambda| \ : \ \lambda \mbox{ is an eigenvalue of
     } B \mbox{ on } V_0\}. $$
\end{lemma}
\begin{pf}
Since $B^T B$ is a real symmetric matrix, its eigenvalues are real. Moreover, for any $\xx \in \C^N$,
we have $(B^T B\xx,\xx)=(B\xx, B\xx) \geq 0$, so these eigenvalues are
non-negative. We have
\begin{equation} \label{eta}
 \eta = \max \{ \ (B^T B \xx, \xx) \, : \,  \xx \in V_0, \ ||\xx||=1 \}
        =  \max \{ \ (B \xx, B\xx) \, : \,  \xx \in V_0, \  ||\xx||=1 \}.
\end{equation}
Now let $\yy \in V_0$ be an eigenvector of $BP$, corresponding to the
eigenvalue $\lambda$, and normalised so that $||\yy||=1$.
Then
$$ |\lambda|^2 = (\lambda \yy, \lambda \yy) = (BP \yy, BP \yy) = (B\zz, B\zz),  $$
where $\zz=P\yy$. But $\zz \in V_0$ since $P$ is column stochastic,
and $||\zz||=1$ since $P$ is orthogonal, so that
$|\lambda|^2 \leq \eta$ as claimed.

Now suppose that $B$ is also a circulant matrix. Let $\yy_j=N^{-1/2}
\vv_j$, where the $\vv_j$ are defined in (\ref{def-vj}). Then the
$\yy_j$ for $1 \leq j \leq N-1$ form an orthonormal basis of
eigenvectors for $B$ on $V_0$. Let $\lambda_j$ be the eigenvalue for
$\yy_j$, and let $k$ be an index such that $|\lambda_k| = \max_{1 \leq
  j \leq N-1} |\lambda_j|$. For any $\xx \in V_0$ with $||\xx||=1$, we
may write $\xx=\sum_{j=1}^{N-1} c_j \yy_j$ with $\sum_{j=1}^{N-1}
|c_j|^2=1$. Then
$$ (B\xx,B\xx)=\sum_{j=1}^{N-1} |c_j|^2 |\lambda_j|^2 \leq
|\lambda_k|^2 =(B \yy_k,B\yy_k), $$
so the maximum in (\ref{eta}) is attained at
$\xx=\yy_k$, giving $\eta=|\lambda_k|^2$.
\end{pf}

\begin{pfThm2-2}
%
For a given $\sigma \in S_N$, we are interested in the eigenvalues
of the matrix $\Phi(1)=m^{-1}B(m,N)Q(\sigma)$, since these are the
eigenvalues of $\Phi(1)$ (where $\Phi(z)$ is the Fredholm matrix of
$\sigma \circ f$) and therefore the isolated eigenvalues in
$\mathrm{Spec}(\mathcal{L}_{f}|_{\operatorname{BV}}).$
By Lemma \ref{shrink-matrix}, it suffices to consider the
eigenvalues of $m^{-1}A(m,N)P(\sigma)$.

The matrix $C(m,N)$ was
obtained from $A(m,N)$ by applying some permutation $\rho$ to the {\em
  rows}. Thus $P(\rho^{-1})A(m,N)=C(m,N)$. For any $\sigma \in S_N$,
the matrix $m^{-1}A(m,N)P(\sigma)=m^{-1}P(\rho)C(m,N)P(\sigma)$ is
conjugate to $m^{-1}C(m,N)P(\sigma \rho^{-1})$, so it suffices to
consider the eigenvalues of the doubly stochastic matrices
$m^{-1}C(m,N)P(\sigma)$ for all $\sigma \in S_N$. We exclude the
eigenvalue 1 associated to the trivial eigenvector $\ee$, so consider
only the eigenvalues on its orthogonal complement $V_0$.

We apply Lemma \ref{circulant-worst} to the doubly stochastic
circulant matrix $B=m^{-1}C(m,N)$, so that $\eta=m^{-1} |\lambda_1 |$
by Proposition \ref{max-circ-eval}. This shows that, for any $\sigma
\in S_N$, each eigenvalue $\lambda$ of $m^{-1}C(m,N)P(\sigma)$
satisfies $|\lambda| \leq m^{-1} \lambda_1$. Thus, in the notation of
Theorem \ref{thm:perm-mix}, we have shown that $\tau_\sigma \leq
\taumax$. Moreover, $\taumax = m^{-1} | \lambda_1| = (-1)^{m-1} m^{-1}
\lambda_1$.

Finally, we must show that each of the values $(-1)^{m-1} e^{2 \pi i
  j/N} \taumax$ and $(-1)^m \taumax$ occurs as an eigenvalue of
$m^{-1} A(m,N)P(\sigma)$ for some $\sigma$. But each of
$m^{-1}C^{(j)}(m,N)$ (for $0 \leq j \leq N-1$) and $m^{-1} C'(m,N)$ is
conjugate to one of these matrices, since $C^{(j)}(m,N)$ and $C'(m,N)$
can be obtained by permuting the rows of $A(m,N)$. In particular,
we have matrices whose eigenvalues include
$ m^{-1} \omega_1^j \lambda_1$ and
$\pm m^{-1} \lambda_1$, as claimed.
\end{pfThm2-2}

We finish this section by noting a further consequence of our discussion of
circulant matrices.

\begin{proposition}
For any $\sigma \in S_N$, the matrix
$A(m,N)P(\sigma)$ has eigenvalue $m$ with (algebraic) multiplicity
1. All its other eigenvalues are algebraic integers of norm $\pm
1$. Composition with $\sigma$ preserves the mixing rate of $f$ (that
is, $\tau_\sigma = 1/m$ in the notation of Section
\ref{sec3:statement}) if and only if these algebraic integers are
roots of unity.
\end{proposition}
\begin{pf}
Clearly the characteristic polynomial of $A(m,N)P(\sigma)$ has integer
coefficients and has leading coefficient $1$, i.e.~its roots are
algebraic integers. If $\lambda$ is any one of these eigenvalues, then
its conjugates are also eigenvalues, and its norm (i.e.~the product of
its conjugates) must be a rational integer. Now the product of the
eigenvalues is $\pm \det(A(m,N)P(\sigma))=\pm \det(C(m,N))\det(P(\rho
\sigma)) = \pm m$ since any permutation matrix has determinant $\pm
1$.  We have the obvious eigenvalue $m$ (with eigenvector $\ee$), so
$m$ has multiplicity 1 as a root of the
characteristic polynomial, and all the other roots must have norm $\pm
1$.

Now if all the eigenvalues $\lambda \neq m$ of $A(m,N)P(\sigma)$ are
roots of unity, we have $|\lambda|=1$. Thus no element of
$\mathrm{Spec}(\mathcal{L}_{f}|_{\operatorname{BV}})$ has modulus
between $m^{-1}$ and $1$, and $\sigma \circ f$ has the same mixing
rate as $f$. Conversely, suppose that $\sigma \circ f$ and $f$ have
the same mixing rate. Then we must have $|\lambda|\leq 1$ for all
eigenvalues $\lambda \neq m$ of $A(m,N)P(\sigma)$. But then all the
conjugates $\lambda'$ of $\lambda$ are again eigenvalues, and hence
satisfy $|\lambda'|\leq 1$. In fact each $|\lambda'|=1$, since the
product of the $\lambda'$ is $\pm 1$. Now any algebraic integer all
of whose conjugates have modulus $1$ must be a root of unity (see
e.g.~\cite[IV, (4,5a)]{FrohlichTaylor}). Hence all the eigenvalues
$\lambda \neq m$ of $A(m,N)P(\sigma)$ are roots of unity.
\end{pf}

\section{Further examples} \label{examples}

We give two examples to demonstrate that the conclusions of Theorems
\ref{new-mixing-theorem} and \ref{thm:perm-mix} do not necessarily
hold if we replace our standard map $f(x)=mx \mod 1$ by other interval maps.
In \S\ref{subshift}, we give an example of a Markov map $f$ where the
proportion of permutations $\sigma \in S_N$ with $\sigma \circ f$ non-mixing
is bounded away from $0$ as $N \to \infty$. Thus the conclusion of
Theorem \ref{new-mixing-theorem} does not hold. In \S\ref{speed-up},
we give a family of interval maps $f$ for which composition with
permutations typically improves the mixing rate, in contrast to
Theorem \ref{thm:perm-mix}(i).

\subsection{An example with many non-mixing permutations} \label{subshift}

Consider the piecewise continuous
function $f : [0,1) \lra [0,1)$ given by
\begin{equation}\label{eq:subshift}
f(x) = \begin{cases}
           2x & \mbox{if } 0 \leq x <\frac{1}{2}, \cr
            x-\frac{1}{2} & \mbox{if } \frac{1}{2} \leq x < 1.
\end{cases}
\end{equation}
Fix $\ell \geq 1$ and divide $[0,1)$ into $N=2\ell$ equal subintervals
$$ I_j = \left[ \frac{j}{2\ell}, \frac{j+1}{2\ell} \right), \quad 0
    \leq j \leq 2\ell-1. $$
For a permutation $\sigma \in S_{2\ell}$ of these subintervals, let
$g=\sigma \circ f$. We have the following result.
\begin{proposition}
The proportion of permutations
$\sigma$ for which $g$ is non-mixing
is bounded away from 0 as $\ell \rightarrow \infty$.
\end{proposition}
\begin{pf}
For any subset $A \subseteq \{0,\ldots, 2\ell-1\}$, define
$\tilde{g}(A) \subseteq \{0,\ldots, 2\ell-1\}$
by
$$ \tilde{g}(A) = \{ \sigma(2j), \sigma(2j+1) \ :\  j \in
    A,\ j<\ell\} \cup \{\sigma(j-\ell) \ : j \in A, j \geq \ell\}.$$
Then, analogously to Proposition \ref{g-union}, we have
$$ g\left( \bigcup_{a \in A} I_a \right) =
      \bigcup_{b \in \tilde{g}(A)} I_b. $$
Note however that Proposition \ref{size-At} no longer holds: for example, if
$A=\{0,\ell,\ell+1\}$ (with $\ell \geq 2$) then
$\tilde{g}(A)=\{0,1\}$ has fewer elements than $A$.

Now if there is some non-empty subset $A$ such that
$\tilde{g}^r(A)\neq \{0,\ldots, 2\ell-1\}$ for all $r \geq 0$
then $g$ is non-mixing.
But if $\sigma$ has the
property that $\sigma(j-\ell)=j$ for some $j \geq \ell$ then,
taking $A=\{j\}$, we have
$\tilde{g}^r(A)=A$ for all $r$. Thus $g$ is 
non-mixing.
We therefore need to investigate the proportion of permutations with
the above property.

Let $1\leq m\leq \ell$ and let $S$ be a subset of $\{\ell,
\ldots,2\ell-1\}$ of size $m$. There are $(2\ell-m)!$ permutations
$\sigma \in S_{2\ell}$ such that $\sigma(j-\ell)=j$ for all $j\in
S$. Moreover, the number of such sets $S$ of size $m$ is ${\ell
  \choose m}$. Thus, by the Inclusion-Exclusion Principle, the
proportion of permutations $\sigma \in S_{2\ell}$ with
$\sigma(j-\ell)=j$ for at least one $j \geq \ell$ is
$$ \sum_{m=1}^\ell (-1)^{m-1} a_m $$
where
$$ a_m = {\ell \choose m} \frac{(2\ell-m)!}{(2\ell)!}. $$
Now the terms in the alternating series are decreasing: for $1 \leq m
< \ell$ we have
$$ \frac{a_{m+1}}{a_m} =
  \frac{ m! (\ell-m)! (2\ell-m-1)! }{ (m+1)! (\ell-m-1)! (2\ell-m)!}
         = \frac{\ell-m}{(m+1)(2\ell-m)} < \frac{1}{2(m+1)}. $$
Hence the required proportion is bounded below by
$$ a_1 - a_2 = \frac{1}{2} - \frac{\ell(\ell-1)}{2(2\ell)(2\ell-1)}
>\frac{3}{8}. $$
So we have proved that, for each $\ell \geq 1$, the function $\sigma
\circ f$ is non-mixing
for more than $3/8$ of
the permutations $\sigma \in S_{2\ell}$.
\end{pf}

We remark that, although the map $f$ in \eqref{eq:subshift} is not
expanding throughout its domain, its second iterate $f^2$ is piecewise
expanding with expansion factor at least 2 everywhere, so our general
discussion of mixing rates can still be applied.

\subsection{An example where permutations speed up mixing}
\label{speed-up}
Consider the following family of
\emph{intermittency maps} $f_{\alpha}:[0,1]\to[0,1]$,
$\alpha\in(0,1)$ given by
\begin{equation}\label{intermittent}
f_{\alpha}(x)=
\begin{cases}
x(1+2^{\alpha}x^{\alpha}) &\textrm{if}\,x\in[0,1/2],\\
2x-1 &\textrm{if}\,x\in(1/2,1].
\end{cases}
\end{equation}
This family has been widely studied \cite{Hu04,LSV99, Young98} and
optimal decay of correlations/speed of convergence to equilibrium
has been established in \cite{Gouezel04}.  In particular, it is
shown that there exists a Banach space $\mathcal{B}$ (e.g., space of
Lipschitz continuous functions), such that
\begin{equation}
||\mathcal{L}_{f}^{n}(\phi)-\rho||_{\mathcal{B}}\leq Cn^{-(1-\frac{1}{\alpha})}||\varphi||_{\mathcal{B}},
\end{equation}
for all $\phi \in \mathcal{B}$ with $\|\phi\|_1=1$.
Moreover this asymptotic in $n$ is optimal within
$\mathcal{B}$. The sub-exponential mixing rate arises since each
$f_{\alpha}$ admits a neutral fixed point at $x=0$, namely
$f'(0)=1$. Thus $f$ is expanding but not uniformly expanding, and the existence of the neutral fixed point inhibits the
mixing. In particular, the functional analytic methods discussed in
Section \ref{backgroud} do not apply since $\lambda=1$ is no
longer an isolated eigenvalue of $\mathcal{L}_f$. i.e. there is no
spectral gap.

We now consider $f_{\alpha}$ composed with a permutation $\sigma\in
S_N$ such that $\sigma \circ f_\alpha$ is topologically mixing. Most
choices of $\sigma$ will not fix the interval $[0,1/N]$, so that
$\sigma \circ f_\alpha$ no longer has a neutral fixed point and is
uniformly expanding on $[0,1]$. Thus $\sigma \circ f_\alpha$ has
bounded variation on $[0,1]$, and it follows from \cite{Keller84, Viana97} that
$\sigma\circ f$ has absolutely continuous invariant measure, with
density in BV. Since the system is uniformly expanding, the operator
$\mathcal{L}_{\sigma\circ f}$ now has a spectral gap, so that rate
of convergence to equilibrium is exponentially fast.

\section*{Acknowledgments}
In preparing this work the authors would like to acknowledge
Congping Lin for numerical assistance. The second-named author
acknowledges support from PREDEX, and the third-named author is
grateful for hospitality from CPT, Marseille and LAMA, University
Paris-East Cr\'{e}teil Val de Marne, Paris.

\end{document}